\documentclass[nohyperref]{article}

\usepackage{microtype}
\usepackage{graphicx}
\usepackage{subfigure}
\usepackage{booktabs} 

\usepackage{hyperref}

\usepackage{url}

\usepackage{pifont}

\usepackage{epsfig}
\usepackage{amssymb}
\usepackage{amsmath}
\usepackage{amsthm}
\usepackage{amsfonts}
\usepackage{bbding}
\usepackage{array}
\usepackage{stfloats,refcount}
\usepackage{paralist}
\usepackage{xargs}                      
\usepackage{caption}

\usepackage{multirow}

\newcolumntype{C}[1]{>{\centering\let\newline\\\arraybackslash\hspace{0pt}}m{#1}}

\newcommand\tagthis{\addtocounter{equation}{1}\tag{\theequation}}

\DeclareMathOperator{\Ocal}{\mathcal{O}}

\newcommand{\E}{\mathbb{E}} 

\newcommand{\sets}[1]{\{#1\}}

\newcommand{\norms}[1]{\Vert#1\Vert}

\newcommand{\iprods}[1]{\langle#1\rangle}

\usepackage{xcolor}

\newcommand{\R}{\mathbb{R}}

\newcommand{\Eproof}{\hfill $\square$}

\newcommand{\zero}[1]{{\boldsymbol{0}}}

\usepackage[accepted]{icml2022}

\usepackage{amsmath}
\usepackage{amssymb}
\usepackage{mathtools}
\usepackage{amsthm}

\usepackage[capitalize,noabbrev]{cleveref}

\theoremstyle{plain}
\newtheorem{theorem}{Theorem}[section]

\newtheorem{lemma}[theorem]{Lemma}
\newtheorem{corollary}[theorem]{Corollary}
\theoremstyle{definition}

\newtheorem{assumption}[theorem]{Assumption}
\theoremstyle{remark}
\newtheorem{remark}[theorem]{Remark}

\usepackage[textsize=tiny]{todonotes}

\icmltitlerunning{Nesterov Accelerated Shuffling Gradient Method for Convex Optimization}

\begin{document}

\twocolumn[
\icmltitle{Nesterov Accelerated Shuffling Gradient Method for Convex Optimization}



\icmlsetsymbol{equal}{*}

\begin{icmlauthorlist}
\icmlauthor{Trang H. Tran}{to1}
\icmlauthor{Katya Scheinberg}{to1}
\icmlauthor{Lam M. Nguyen}{to2}
\end{icmlauthorlist}

\icmlaffiliation{to1}{School of Operations Research and Information Engineering, Cornell University, Ithaca, NY, USA.}
\icmlaffiliation{to2}{IBM Research, Thomas J. Watson Research Center, Yorktown Heights, NY, USA}

\icmlcorrespondingauthor{Lam M. Nguyen}{LamNguyen.MLTD@ibm.com}

\icmlkeywords{Machine Learning, ICML}

\vskip 0.3in
]



\printAffiliationsAndNotice{}  

\begin{abstract}
In this paper, we propose Nesterov Accelerated Shuffling Gradient (NASG), a new algorithm for the convex finite-sum minimization problems. Our method integrates the traditional Nesterov's acceleration momentum with different shuffling sampling schemes. We show that our algorithm has an improved rate of $\Ocal(1/T)$ using unified shuffling schemes, where $T$ is the number of epochs. This rate is better than that of any other shuffling gradient methods in convex regime. Our convergence analysis does not require an assumption on bounded domain or a bounded gradient condition. For randomized shuffling schemes, we improve the convergence bound further. When employing some initial condition, we show that our method converges faster near the small neighborhood of the solution. Numerical simulations demonstrate the efficiency of our algorithm. 
\end{abstract}

\section{Introduction}\label{sec_intro}

We consider the finite-sum optimization problem: 
\begin{equation}\label{ERM_problem_01}
    \min_{w \in \mathbb{R}^d} \left\{ F(w) := \frac{1}{n} \sum_{i=1}^n f(w ; i)  \right\},
\end{equation}
where the objective function $F: \R^d \to \R$ is smooth and convex, and each individual functions $f_i$ is smooth. This standard problem arises in most machine learning tasks, including logistic regression, multi-kernel learning, and some neural networks. The major challenge in solving \eqref{ERM_problem_01} often comes from the high dimension space and a large number of components $n$. Therefore, deterministic methods relying on full gradients are usually inefficient to solve this problem \citep{sra2012optimization, Bottou2018}. 

\textbf{\textit{SGD and Shuffling SGD.}}
Stochastic Gradient Descent (SGD) \cite{RM1951} and its stochastic first-order variants have been widely used to solve \eqref{ERM_problem_01} thanks to its scalability and efficiency in dealing with large-scale problems \citep{AdaGrad,Kingma2014,Bottou2018,Nguyen2018_sgdhogwild}. 

At each iteration SGD samples an index $i$ uniformly from the set $\{1, \dots, n\}$, and uses the stochastic gradient $\nabla f_i$ to update the weight.
While the uniformly independent sampling of $i$ plays an important role in our theoretical understanding of SGD, practical heuristics often use without-replacement sampling schemes (also known as shuffling sampling schemes). 
These methods depend on some random or deterministic permutations of the index set  $\{1, 2, \dots, n\}$ and apply incremental gradient updates using these permutation order. A collection of such $n$ individual updates is called an epoch, or a pass over all the data. 
The most popular method in this class is Random Reshuffling, which creates a new random permutation at the beginning of each epoch. Other important methods  include Single Shuffling (which uses the same (random) permutation for each epoch) and Incremental Gradient (which uses a deterministic order of indices). 
In this paper, the term Shuffling SGD refers to SGD method using \textit{any} data permutations, which includes the three special schemes described above. 

Empirical studies show that shuffling sampling schemes usually provide a faster convergence than SGD \citep{bottou2009curiously}. However, due to the lack of statistical independence, analyzing these shuffling variants is often more challenging than the identically distributed version. Recent works have shown theoretical improvement for shuffling schemes over  SGD in terms of the number of epochs needed to converge to an $\epsilon$-accurate solution\footnote{We define an $\epsilon$-accurate solution as a point $x \in \R^d$ that satisfies $F(x) - F(x_*) \leq \epsilon$ for convex settings (where $x_*$ is a minimizer of $F$ and the statement may hold in expectation). }  \citep{Gurbuzbalaban2019,haochen2019random,safran2020good,nagaraj2019sgd,rajput2020closing,nguyen2020unified,mishchenko2020random, ahn2020sgd}.
In particular, in a general convex setting, shuffling sampling schemes improve the convergence rate of SGD from $\Ocal(1/\sqrt{T})$ to $\Ocal(1/T^{2/3})$ 
in terms of the number of effective data passes $T$
\citep{nguyen2020unified,mishchenko2020random}. Thanks to their theoretical and empirical advantage, Random Reshuffling and its  variants are becoming the methods of choice for practical implementation of machine learning optimization algorithms.

\textbf{\textit{Nesterov's Accelerated Gradient (NAG).}} On the other hand, one of the most beautiful idea in convex optimization is the Nesterov's accelerated momentum technique, which was originally proposed in \citep{Nesterov1983}. The method, shown in Algorithm \ref{alg_nesterov} for deterministic setting, achieves a much better convergence rate of $\Ocal(1/T^2)$ than the convergence rate of Gradient Descent $\Ocal(1/T)$ in convex regimes, where $
T$ is the total number of iterations. Note that the application of deterministic NAG requires a full gradient computation, i.e $n$ component gradients in each iteration, this $T$ is the same as the number of epochs.  

\begin{algorithm}[hpt!]
   \caption{Nesterov's Accelerated Gradient (NAG)}\label{alg_nesterov}
\begin{algorithmic}[1]
   \STATE {\bfseries Initialization:} Choose an initial point $x_0, y_0\in \mathbb{R}^d$.
   \FOR{$t=1,2,\cdots,T$}
   \STATE Let $x^{(t)} := y^{(t-1)} - \alpha^{(t)} \nabla F ( y^{(t-1)} ) $
   \STATE Compute $y^{(t)} := x^{(t)} + \frac{t-1}{t+2} ( x^{(t)} - x^{(t-1)}  )$
   \ENDFOR
\end{algorithmic}
\end{algorithm} 

In the last two decades, researchers have made efforts to leverage this acceleration technique to the stochastic settings. 
It is well known that Stochastic Gradient Descent has the convergence rate of $\Ocal(1/\sqrt{K})$  where $K$ is the number of iterations\footnote{To make fair comparisons, we use $K$ for the  iteration of SGD. Note that $K$ is the number of individual gradient evaluations, and it is equivalent to $nT$ in other methods that use $T$ data passes.}.
\citet{pmlr-v89-vaswani19a} proposes to use a new assumption called the Strong Growth Condition for which they can prove an accelerated rate of SGD with Nesterov's momentum. This 
condition implies that the stochastic gradients (and, hence, its variance) converge to zero at the optimum \cite{schmidt2013fast, pmlr-v89-vaswani19a}.
However, without a strong assumption on the gradient oracle (i.e. without assuming that the variance goes to zero), no work has been able to prove a better convergence rate for SGD with Nesterov's momentum over the ordinary results of SGD \citep{NIPS_hu2009,Lan2012}. This background along with the theoretical advances of Shuffling SGD motivates the central question of our paper:
\begin{center}
\textit{Can we use Nesterov's momentum technique for Shuffling SGD to improve the convergence rate using only standard assumptions (e.g. without assuming vanishing variance)?}
\end{center}

We answer this question positively in this paper; our results are summarized below.



\textbf{{Summary of our contributions}}. 
\begin{itemize}
    \item We propose Nesterov Accelerated Shuffling Gradient (NASG) method, a new algorithm to approximate the solution of the convex minimization problem \eqref{ERM_problem_01}. Our method integrates the well-known Nesterov's acceleration technique with shuffling sampling strategies. In stead of the traditional practice that add momentum term in each iteration, we adopt a new approach that integrates the momentum for each training epoch. 
    \item We establish the convergence analysis for our algorithm in the convex setting using standard assumptions, i.e. generalized bounded variance or convex component functions. 
    Our method achieves an improved rate of $\Ocal(1/T)$ in terms of the number of epochs for the unified shuffling schemes.  
    We also investigate the randomized schemes (including Random Reshuffling and Single Shuffling) and improve a factor of $n$ in the convergence bound.
    Moreover, our convergence results work for the last iterate returned by the algorithm, which is more practical than previous works for the average iterate. 
    \item We test our algorithms via numerical simulations on various machine learning tasks and compare them with other stochastic first order methods. Our tests have shown good overall performance of the new algorithms.
\end{itemize}

\textbf{Related work.}
Let us briefly review the most related works to our methods studied in this paper.

\textbf{\textit{Shuffling SGD schemes.}} In the big data machine learning setting, Random Reshuffling and Single Shuffling are more favorable than plain SGD thanks to their better practical performance and simple implementation \citep{bottou2009curiously,bottou2012stochastic, Recht2011}. While the convergence properties of SGD are well-understood in literature, the theoretical analysis for the randomized shuffling schemes remained challenging for a long period of time. A natural reason behind this problem is the lack of conditionally unbiased gradients: $\E \left [\nabla f(y_i^{(t)}; \pi_i^{(t)}) \right] \neq \nabla F(y_i^{(t)}) ,$ where $t$ is the current epoch. 
Recently, researchers have made progress in the analysis of convergence
rates of randomized shuffling techniques \citep{Gurbuzbalaban2019,haochen2019random,safran2020good,nagaraj2019sgd, ahn2020sgd}. 
with the majority of these works devoted to the strongly convex case (with a bounded gradient or bounded domain assumption). The best known convergence rate in this case is  $\Ocal(1/(nT)^2 + 1/(nT^3))$ where $T$ is the number of epochs. This result matches the lower bound rate in  \citep{safran2020good} up to some constant factor.

\renewcommand{\arraystretch}{1.5}
\begin{table*}
\caption{Number of individual gradient evaluations needed by SGD-type algorithms to reach an $\epsilon$-accurate solution $x$ that satisfies $F(x) - F(x_*) \leq \epsilon$. 
In this table, $L$ is the Lipschitz constant in Assumption \ref{ass_smooth}, $\sigma_*^2$ is the variance at the minimizer $x_*$ defined in \eqref{defn_finite}. Finally, $\Delta := \| \tilde{x}_0 -  x_{*} \|^2$ is the squared distance from the initial point $\tilde{x}_0$ to the minimizer $x_*$.
}\label{table_1}
\begin{center}
\begin{tabular}{ |l|l|l| } 
 \hline
 Algorithms & Complexity & References \\ 
 \hline
 Standard SGD\textcolor{red}{$^{(1)}$} &{$\Ocal\left(\frac{\Delta_0^2 + G^2}{\epsilon^2}\right)$} \textcolor{red}{$^{(1)}$}  &
 \citep{Nemirovski2009,pmlr-v28-shamir13}  \\ 
 \hline
 SGD - Unified Schemes\textcolor{red}{$^{(2)}$} &{$\Ocal\left(\frac{nL\Delta}{\epsilon} + \frac{n\sqrt{L}\sigma_* \Delta}{\epsilon^{3/2}}\right)$ } & 
 \citep{mishchenko2020random,nguyen2020unified} \\ 
 \hline
 SGD - Randomized Schemes\textcolor{red}{$^{(3)}$} &{$\Ocal\left(\frac{nL\Delta}{\epsilon} + \frac{\sqrt{nL}\sigma_* \Delta }{\epsilon^{3/2}}\right)$ }  & 
 \citep{mishchenko2020random}\\
 \hline
 
\textcolor{blue}{NASG - Unified Schemes}  &{$\Ocal\left( \frac{nL\Delta}{\epsilon} + \frac{n\sigma_*^2}{L\epsilon} \right) $ } & 
 \textcolor{blue}{(This work, Theorem \ref{thm_convex_1} and Corollary \ref{cor_comp_unified}) } \\ 
 \hline

 \textcolor{blue}{NASG - Randomized Schemes}\textcolor{red}{$^{(3)}$}  &{$\Ocal\left( \frac{nL\Delta}{\epsilon} + \frac{\sigma_*^2}{L\epsilon} \right)$ }& 
\textcolor{blue}{(This work, Theorem \ref{thm_convex_RR} and Corollary \ref{cor_comp_RR})} \\ 
 \hline
\end{tabular}
\end{center}
\footnotesize{\textcolor{red}{$^{(1)}$} Standard results for SGD in literature often use a different set of assumptions  from the one in this paper (e.g. bounded domain that $\| x -  x_{*} \|^2 \leq \Delta_0$ for each iterate $x$ and/or  bounded gradient that $\E[\|\nabla f(x;i)\|] \leq G^2 $). We report the corresponding complexity for a rough comparison. 
\textcolor{red}{$^{(2)}$} \citep{mishchenko2020random} shows a bound for Incremental Gradient while  \citep{nguyen2020unified} has a unified setting. We translate these results for Unified Schemes from these references to our convex setting. 
\textcolor{red}{$^{(3)}$} While using the same set of assumptions, the convergence criteria for randomized schemes is in expectation form: $\E[F(x) - F(x_*)] \leq \epsilon$.}
\end{table*}

In the convex regime, most dominant results are originally derived for the deterministic Incremental Gradient scheme \citep{nedic2001convergence,nedic2001incremental}. More recent works investigate convergence theory for various shuffling schemes  \citep{shamir2016without,mishchenko2020random,nguyen2020unified}, where \citet{nguyen2020unified} provides a unified approach to different shuffling schemes and  proves the convergence rate of $\Ocal(1/ T^{2/3})$. When a randomized scheme is applied (Random Reshuffling or Single Shuffling), the bound in expectation improves to $\Ocal(1/T + 1/(n^{1/3} T^{2/3}))$. 
For a comparison, our Algorithm~\ref{shuffling_nesterov_02} developed in this paper 
achieves a deterministic convergence rate of $\Ocal(1/T)$ for the same setting under standard assumptions. 
The computational complexity for these methods are in Table \ref{table_1}.

In the meantime, a popular line of research involves variance reduction technique, which have shown encouraging performance for machine learning (e.g., SAG \citep{SAG}, SAGA \citep{SAGA}, SVRG \citep{SVRG} and SARAH \citep{Nguyen2017sarah}). These methods need to either compute or store a full gradient or a large batch of gradient. This  plays an important role in reducing the variance and therefore, is the key factor for these methods. 
However, the update of SGD, Shuffling SGD and our Algorithm \ref{shuffling_nesterov_02} does not require full gradient evaluation at any stage.
Thus, our new Algorithm \ref{shuffling_nesterov_02} belongs to the class of Shuffling SGD which deviates from variance reduction methods. 

\textbf{\textit{Momentum Techniques.}}
The most popular and successful momentum techniques include the classical Heavy-ball method \citep{Polyak1964} and Nesterov’s acceleration gradient (NAG) \citep{Nesterov1983, Nesterov2004}. Although these two methods are different, they both receive great attention in the optimization community \citep{NIPS_hu2009,Lan2012,pmlr-v28-sutskever13,JMLR_yuan2016,dozat2016incorporating}.
Nesterov’s acceleration method is well-known for its improved convergence rate of $\Ocal(1/T^2)$  (versus the $\Ocal(1/T)$ of Gradient Descent) for general smooth convex functions in the deterministic setting, where $T$ is the number of iterations.

On the other hand, \citet{springer_devolder2014} and  \citet{siam_lessard2016} suggest that Nesterov’s acceleration is not robust to the errors in gradient and its performance may be worse than gradient descent due to error accumulation. A more recent work \cite{corr_liu2018} argues that stochastic NAG does not provide acceleration over ordinary SGD in general, and may diverge for step sizes that guarantee convergence of SGD. These observations further motivate our algorithmic design for the Shuffling SGD with Nesterov's momentum in Section \ref{sec:shuffling_nesterov}.

\section{Nesterov Accelerated Shuffling Gradient Method}\label{sec:shuffling_nesterov}
In this section, we describe our new shuffling gradient algorithm with Nesterov's momentum in Algorithm~\ref{shuffling_nesterov_02}. 

Before we start, it should be noted that the classical approach in stochastic NAG literature is applying the momentum term for each iteration \citep{NIPS_hu2009,Lan2012, pmlr-v33-zhong14, pmlr-v89-vaswani19a}. However, empirical evidence have shown that Nesterov’s acceleration may not be superior when inexact gradients are used, and the reason might be error accumulation \cite{springer_devolder2014,corr_liu2018}.
In addition, while SGD has access to an unbiased estimator for the full gradient, shuffling gradient schemes generally do not have this property. In consequence, updating the momentum at each inner iteration is less preferable since it could make the estimator deviate from the true gradient and further accumulate errors.

Based on these observations, we adopt a different approach to update the Nesterov's momentum after each epoch which consists of $n$ gradients. This practice allows our method to approximate the full gradient more accurately while still maintains the effectiveness of the momentum technique. It is also consistent with the application of Heavy-ball method and proximal operator for shuffling schemes in recent literature \cite{pmlr-v139-tran21b, mishchenko2021proximal}.
Our algorithm is presented below.

\begin{algorithm}[hpt!]
   \caption{Nesterov Accelerated Shuffling Gradient (NASG) Method}\label{shuffling_nesterov_02}
\begin{algorithmic}[1]
   \STATE {\bfseries Initialization:} Choose an initial point $\tilde{x}_0, \tilde{y}_0 \in \mathbb{R}^d$.
   \FOR{$t=1,2,\cdots,T $}
   \STATE Set $y_0^{(t)} := \tilde{y}_{t-1}$;
   \STATE Generate any permutation  $\pi^{(t)}$ of $[n]$ (either deterministic or random);
   \FOR{$i = 1,\cdots, n$} 
    \STATE Update $y_{i}^{(t)} := y_{i-1}^{(t)} - \eta_i^{(t)} \nabla f ( y_{i-1}^{(t)} ; \pi^{(t)} ( i ) )$; 
   \ENDFOR 
   \STATE Set $\tilde{x}_t := y_{n}^{(t)}$;
   \STATE Update $\tilde{y}_{t} := \tilde{x}_{t} + \gamma_t ( \tilde{x}_{t} - \tilde{x}_{t-1} )$; 
   \ENDFOR 
\end{algorithmic}
\end{algorithm} 
\textbf{\textit{Algorithm Description.}} 
In each epoch $t$, our method first performs $n$ consecutive individual gradient updates in variable $y_i^{(t)}$ following a permutation $\pi^{(t)}$ of the index set $\{ 1, \dots, n\}$. At the end of each epoch, it applied the Nesterov's momentum update using an auxiliary variable $\tilde{x}_t$. 
The choice of learning rate $\eta_i^{(t)}$ is further specified in our theoretical analysis.  

The per-iteration complexity of Algorithm \ref{shuffling_nesterov_02} is the same as standard shuffling gradient schemes  \citep{shamir2016without}. In addition, our algorithm only requires a storage cost of $\Ocal(d)$, which is similar to that of standard SGD. 
Note that the implementation of our method  requires neither full gradient computation nor a large batch of gradient computation at any point.
Our convergence guarantee in Theorem \ref{thm_convex_1} and Theorem \ref{thm_convex_2} for unified shuffling scheme holds for any permutation of $\{1,2,\cdots, n\}$, including deterministic and random ones. 
Therefore, our method works for any shuffling strategy, including Incremental Gradient, Single Shuffling, and Random Reshuffling. 

\textbf{\textit{Comparison with Nesterov’s Accelerated Gradient.}} Let us recall that deterministic NAG has an update of full gradient computation from $y^{(t-1)}$ to $x^{(t)}$. We can write this update in a different way, where $y^{(t-1)} = y_{0}^{(t)}$ and each component gradient at $y_{0}^{(t)}$ is gradually computed and subtracted from the starting point: 

\begin{algorithm}[hpt!]
\begin{algorithmic}[1]
    \setcounterref{ALC@line}{defn_finite}
   \FOR{$i = 1,\cdots, n$}
    \STATE Update $y_{i}^{(t)} := y_{i-1}^{(t)} - \eta_i^{(t)} \nabla f ( \textcolor{blue}{y_{0}^{(t)}} ; \pi^{(t)} ( i ) )$; 
   \ENDFOR
\end{algorithmic}
\end{algorithm} 
With an appropriate choice of learning rates, at the end of an epoch, 
the output $y_{n}^{(t)}$ in this representation is identical to the output $x^{(t)}$ of deterministic NAG algorithm. This illustrates  the comparison between traditional NAG and our method. While  NAG only update the weights after a full gradient computation, our method gradually updates and makes movement after each component evaluations. 

\begin{figure}[H] 
\centering
\includegraphics[height=0.15\textheight]{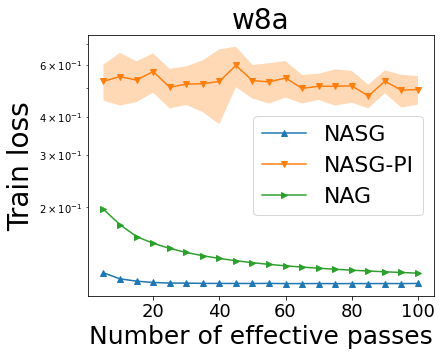}
\includegraphics[height=0.15\textheight]{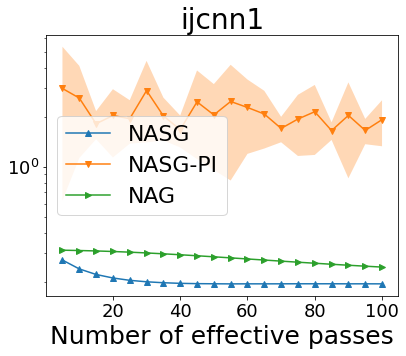}
\caption{Comparisons of the training loss for \texttt{w8a} and \texttt{ijcnn1} datasets. NAG denotes the deterministic Nesterov's Accelerated  Gradient. NASG is our method, while NASG-PI is the stochastic shuffling version that applies Nesterov's momentum each \textit{inner} iteration. We apply random reshuffling schemes for the stochastic algorithms.}
\label{fig_demo}
\end{figure}

In order to motivate our Algorithm \ref{shuffling_nesterov_02}, we conduct a small binary classification experiment and demonstrate the behaviour of NAG and the stochastic momentum methods. The details of  the settings are delayed to Section \ref{subsec:exp_binary}. Figure \ref{fig_demo} shows that applying Nesterov's momentum term for each iteration may accumulate errors and lead to a poor result. While the deterministic NAG converges and slowly decreases the loss, our stochastic version works faster and achieves an overall better performance, when the number of data $n$ is large.  

\section{Technical Settings}\label{sec:technical_settings}
\subsection{Theoretical Assumptions}\label{subsec:technical_assumptions}
We analyze the convergence of Algorithm~\ref{shuffling_nesterov_02} under standard assumptions, which are presented below. Our first assumption is that 
all component functions are  $L$-smooth.
\begin{assumption}
\label{ass_smooth}
We assume that $f(\cdot;i)$ is $L$-smooth for every $i \in [n]$, i.e., there exists a constant $L > 0$ such that, $\forall x,y \in \mathbb{R}^d$, 
\begin{align*}
\| \nabla f(x;i) - \nabla f(y;i) \| \leq L \| x - y \|, \ i \in [n].  \tagthis\label{eq:Lsmooth_basic}
\end{align*}
\end{assumption}


Assumption \ref{ass_smooth} implies that $F$ is also $L$-smooth. This assumption is widely used in literature for gradient-type methods in both stochastic and
deterministic settings. Then, by a well-known property of $L$-smooth functions in \cite{Nesterov2004}, we have, $\forall x, y \in \mathbb{R}^d$, 
\begin{align*}
F(x) &\leq F(y) + \langle \nabla F(y),(x-y) \rangle + \frac{L}{2}\|x-y\|^2. \tagthis\label{eq:Lsmooth}
\end{align*}

Our main results for Algorithm~\ref{shuffling_nesterov_02} is in convex setting which requires the following assumption.
\begin{assumption}\label{ass_convex}
$f(\cdot;i)$ is convex for every $i \in [n]$, i.e., $\forall x, y \in \mathbb{R}^d$, 
\begin{gather*}
f(x;i)  - f(y;i) \geq \langle \nabla f(y;i),(x - y) \rangle, \ i \in [n]. 
\end{gather*}
\end{assumption}



When $F$ is convex, we also assume that the existence of a minimizer for $F$. 
Note that $F$ could have more than one minimizer. Therefore, in this paper, we let $x_*$ be any minimizer of $F$ and consider 
the corresponding variance of $F$ at $x_{*}$:
\begin{equation}\label{defn_finite}
\sigma_{*}^2 := \frac{1}{n}\sum_{i=1}^{n}\Vert \nabla{f}(x_{*}; i) \Vert^2 \in [0, +\infty).
\end{equation}



Alternatively, when each individual function $f(\cdot;i)$ is not necessary convex and $F$ is convex, we consider the following assumption:
\begin{assumption}\label{ass_bounded_variance}
\textbf{(Generalized bounded variance)} 
There exist two non-negative and finite constants $\Theta$ and $\sigma$ such that for any $x \in \mathbb{R}^d$ we have
\begin{equation}\label{eq:bounded_variance}
 \frac{1}{n} \sum_{i=1}^n \norms{\nabla f(x; i) - \nabla F(x)}^2 \leq \Theta \| \nabla F(x) \|^2 + \sigma^2. 
\end{equation}
\end{assumption}
Assumption \ref{ass_bounded_variance} reduces to the standard bounded variance condition if $\Theta = 0$. 
Therefore, it is more general than the bounded variance assumption, which is often used in stochastic optimization \cite{Bottou2018}.

\subsection{Basic Derivations}\label{subsec:basic_derivations}
In this section, we provide some key derivations for Algorithm \ref{shuffling_nesterov_02}. 
From the update of our algorithm and the choice $\gamma_t = \frac{t-1}{t+2}$, we have for $t \geq 1$
\begin{align*}
    \tilde{y}_{t} := \tilde{x}_{t} + \frac{t-1}{t+2} ( \tilde{x}_{t} - \tilde{x}_{t-1} ). \tagthis \label{eq_005}
\end{align*}
Similar to the original Nesterov's momentum technique, we use two following auxiliary variables in our analysis:
\begin{align*}
    \theta^{(t)} = \frac{2}{t+2} \in (0,1),  \text{ and } v^{(t)} = \frac{t+1}{2} \tilde{x}_{t} - \frac{t-1}{2} \tilde{x}_{t-1},
\end{align*}
for $t \geq 1$. We also use the convention that $\theta^{(0)} = 1$ and $v^{(0)} = \tilde{x}_{0}$.
This is equivalent to
\begin{align*}
    \tilde{x}_{t} &= \frac{2}{t+1} v^{(t)} + \frac{t-1}{t+1} \tilde{x}_{t-1} \\
    &= \theta^{(t-1)} v^{(t)} + (1 - \theta^{(t-1)}) \tilde{x}_{t-1}. \tagthis \label{eq_006_0}
\end{align*}
This property shows that $\tilde{x}_{t}$ is a convex combination of $v^{(t)}$ and $\tilde{x}_{t-1}$ for every iteration $t \geq 1$. Using equation \eqref{eq_005}, we further have 
\begin{align*}
    \tilde{y}_{t} & = \tilde{x}_{t} + \frac{t-1}{t+2} ( \tilde{x}_{t} - \tilde{x}_{t-1} ) \\
    &= \frac{t+1}{t+2} \tilde{x}_{t} - \frac{t-1}{t+2} \tilde{x}_{t-1} + \frac{t}{t+2} \tilde{x}_{t} \\
    &= \frac{2}{t+2} \left(  \frac{t+1}{2} \tilde{x}_{t} - \frac{t-1}{2} \tilde{x}_{t-1} \right) + \left(1 - \frac{2}{t+2}  \right) \tilde{x}_{t} \\
    &= \theta^{(t)} v^{(t)} + (1 - \theta^{(t)}) \tilde{x}_{t}. \tagthis \label{eq_006}
\end{align*}
Again, $\tilde{y}_{t}$ is a convex combination of $v^{(t)}$ and $\tilde{x}_{t}$, however with a slightly different parameter $\theta^{(t)}$ instead of $\theta^{(t-1)}$. 

These key derivations play an important role in the theoretical analysis of Algorithm \ref{shuffling_nesterov_02}. They help explain why Nesterov's momentum can achieve a better convergence rate when the objective function $F$ is convex. Indeed, using convexity of $F$ we have the following property: 
\begin{align*}
     &F(y) + \langle \nabla F ( y ) ,  (1 - \theta) x + \theta x_{*} - y \rangle\\
     \leq \ &F ( (1 - \theta) x + \theta x_{*} ) \leq (1 - \theta) F ( x ) + \theta F ( x_{*} )
\end{align*}
for any $x \in \mathbb{R}^d$, $y \in \mathbb{R}^d$, and $\theta \in [0,1]$. The application of this inquality is the central idea behind our theoretical results, which are presented in the next section. 

\section{Theoretical Analysis}
\subsection{Convergence Rate for Unified Shuffling Scheme}

In this section, we investigate the theoretical performance of Algorithm~\ref{shuffling_nesterov_02} using unified  shuffling strategy, i.e. using an arbitrary permutation $\pi^{(t)}$ in any of the epoch $t=1, 2, \dots, T$. These permutations can be random or deterministic, however our results hold deterministically regardless of the choice of permutation. 
We first establish the convergence for Algorithm~\ref{shuffling_nesterov_02} under the condition that {all the component functions are convex}.

\begin{theorem}[Convex component functions]\label{thm_convex_1}
Suppose that Assumption \ref{ass_smooth} and \ref{ass_convex}  hold for \eqref{ERM_problem_01}. 
Let $\sets{x_i^{(t)}}$ be generated by  Algorithm~\ref{shuffling_nesterov_02} with parameter $\gamma_t = \frac{t-1}{t+2}$, the learning rate $\eta_i^{(t)} := \frac{\eta_t}{n} > 0$ for $\eta_t = \frac{k\alpha^t}{LT} \leq \frac{1}{L}$ where $k = \frac{1}{e \alpha \sqrt[3]{12}} > 0$ and $\alpha = 1+ \frac{1}{T}>0$.
Then for $T \geq 2$ we have
\begin{align}
    F ( \tilde{x}_{T} ) - F ( x_{*} ) 
    &\leq \frac{4 \sigma_*^2 }{9  LT}   
    + \frac{2Le \sqrt[3]{12}}{T}  \| \tilde{x}_0 -  x_{*} \|^2. \label{eq_thm_convex_1}
\end{align}

\end{theorem}
\begin{remark}\label{re:main_thm}
The convergence rate of Algorithm~\ref{shuffling_nesterov_02} is exactly expressed as
\begin{equation*}
\Ocal\left(\frac{ \sigma_*^2/L + L\| \tilde{x}_0 -  x_{*} \|^2  }{T} \right),
\end{equation*}
which is better than the state-of-the-art rate in the literature \citep{mishchenko2020random,nguyen2020unified} in term of $T$ for convex problems with general shuffling-type strategies. Translating this convergence rate to computational complexity, we get the results in Table \ref{table_1}. We provide the proof of Theorem \ref{thm_convex_1} and its complexity in the Appendix.
\end{remark}
When the component functions are not necessarily convex, we establish the convergence for Algorithm~\ref{shuffling_nesterov_02} under the { convexity of $F$  and the generalized bounded variance assumption}.

\begin{theorem}[Generalized Bounded Variance]\label{thm_convex_2}
Suppose that Assumption \ref{ass_smooth} and \ref{ass_bounded_variance} hold for \eqref{ERM_problem_01}. In addition, we assume that $F$ is convex.
Let $\sets{x_i^{(t)}}$ be generated by  Algorithm~\ref{shuffling_nesterov_02} with parameter $\gamma_t = \frac{t-1}{t+2}$, the learning rate $\eta_i^{(t)} := \frac{\eta_t}{n} > 0$ for $\eta_t = \frac{k\alpha^t}{LT} \leq \frac{1}{L}$ where $k = \frac{1}{e \alpha \sqrt[3]{2(6\Theta + 7)}} > 0$ and $\alpha = 1+ \frac{1}{T} >0$.
Then for $T \geq 2$, $ F ( \tilde{x}_{T} ) - F ( x_{*} ) $ is upper bounded by   
\begin{align*}
    \frac{8\sigma^2 }{3(6\Theta + 7)LT} 
    + \frac{2Le  \sqrt[3]{2(6\Theta + 7)} }{T}  \| \tilde{x}_0 -  x_{*} \|^2 . 
\end{align*}
\end{theorem}
The convergence rate of Theorem \ref{thm_convex_2} is expressed as
\begin{equation*}
\Ocal\left(\frac{ \sigma^2/(\Theta L) + L \Theta^{1/3} \| \tilde{x}_0 -  x_{*} \|^2  }{T} \right),
\end{equation*}
which is similar to the convergence rate $\Ocal\left(1/T\right)$ of Theorem \ref{thm_convex_1}. We defer the proof of Theorem \ref{thm_convex_RR} to Appendix.

\begin{remark}[Convergence guarantee]
Our convergence bounds in Theorem \ref{thm_convex_1} and \ref{thm_convex_2} hold in a deterministic sense.
This convergence criteria for Algorithm~\ref{shuffling_nesterov_02} is significantly stronger than the standard criteria in expectation for other SGD-type algorithm in literature recently \citep{ghadimi2013stochastic, pmlr-v28-shamir13}. This improvement is made thanks to the unique structure of the Nesterov's acceleration applied to Shuffling schemes in our Algorithm \ref{shuffling_nesterov_02}. In addition, our results hold for the last iterate $x_T$, which 
matches the practical heuristics more than previous results that hold for an average $\tilde{x}$ of training weights $x_1, \dots, x_T$ \citep{Polyak1992, ghadimi2013stochastic}.
\end{remark}

\subsection{Convergence Rate for Randomized Schemes}
We continue to present the theoretical result of Algorithm \ref{shuffling_nesterov_02} specifically for Randomized Schemes, namely Random Reshuffling and Single Shuffling schemes where random permutation(s) are generated for the update of Algorithm \ref{shuffling_nesterov_02}. 
Our next Theorem \ref{thm_convex_RR} uses the assumption that all the component functions are convex.


\begin{theorem}[Randomized Schemes]\label{thm_convex_RR}
Suppose that Assumption \ref{ass_smooth} and \ref{ass_convex} hold for \eqref{ERM_problem_01}. 
Let $\sets{x_i^{(t)}}$ be generated by  Algorithm~\ref{shuffling_nesterov_02} under a \textbf{randomized scheme} with parameter $\gamma_t = \frac{t-1}{t+2}$, the learning rate $\eta_i^{(t)} := \frac{\eta_t}{n} > 0$ for $\eta_t = \frac{k\alpha^t}{LT} \leq \frac{1}{L}$ where $k = \frac{1}{e \alpha \sqrt[3]{12}} > 0$ and $\alpha = 1+ \frac{1}{T}>0$.
Then for $T \geq 2$, 
we have
\begin{align}
    \E[F ( \tilde{x}_{T} ) - F ( x_{*} ) ]    &\leq 
    \frac{8 \sigma_*^2 }{27n LT}   
    + \frac{2Le \sqrt[3]{12}}{T}  \| \tilde{x}_0 -  x_{*} \|^2.
\end{align}
\end{theorem}

\begin{remark}[Randomized Schemes]
The convergence rate of Theorem \ref{thm_convex_RR} is expressed as
\begin{equation*}
\Ocal\left(\frac{ \sigma_*^2/L  }{nT} + \frac{L \| \tilde{x}_0 -  x_{*} \|^2  }{T} \right),
\end{equation*}
which 
is better than the state-of-the-art rate for randomized schemes in the literature \citep{mishchenko2020random,nguyen2020unified} for convex problems.

Comparing to the unified case, our result allows a reduction in the first term of the bound by a factor of $n$. This fact is essentially helpful in machine learning applications where the number of data $n$ is large. Furthermore, in practice randomized schemes offer a lot of improvements when the variance at the optimizer $\sigma_*^2$ can be large. 
Similar to the previous theorems, our result in Theorem \ref{thm_convex_RR} holds for the last iterate $x_T$, which matches the practical heuristics. We defer the proof of this theorem to the Appendix.
\end{remark}


\subsection{Improved Convergence Rate with Initial Condition}

In this section, we consider an initial condition where the iterate of our algorithm is in a small neighborhood of the optimal point. Let us note that the minimizer of $F$ may not be unique, hence we only requires this assumption for some minimizer $x_*$. 
\begin{remark}\label{rem_unified}
Let us assume that $\| \tilde{x}_0 -  x_{*} \| \leq \frac{E}{\sqrt{n}}$ where $\tilde{x}_0$ be the initial point and $E>0$ be a constant. For the same conditions as in Theorem~\ref{thm_convex_1}, i.e. component convexity, 
we have
\begin{align}
     F ( \tilde{x}_{T} ) - F ( x_{*} ) 
    \leq \Ocal\left(\frac{ \sigma_*^2/L + LE^2  }{n^{3/4}T} \right), 
    \label{eq_rem_unified}
\end{align} 
which has an improvement of $n^{3/4}$ over the plain setting of Theorem \ref{thm_convex_1}. This fact suggests that the algorithm may converge faster when it reaches a small neighborhood of the solution set. 
The proof of this Remark requires some modifications from Theorem~\ref{thm_convex_1}, and is presented in Appendix.  
\end{remark}
\begin{figure*}[b]
\centering
\includegraphics[width=0.33\textwidth]{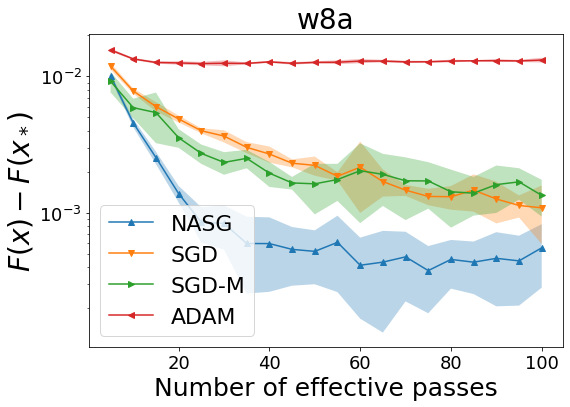}
\includegraphics[width=0.33\textwidth]{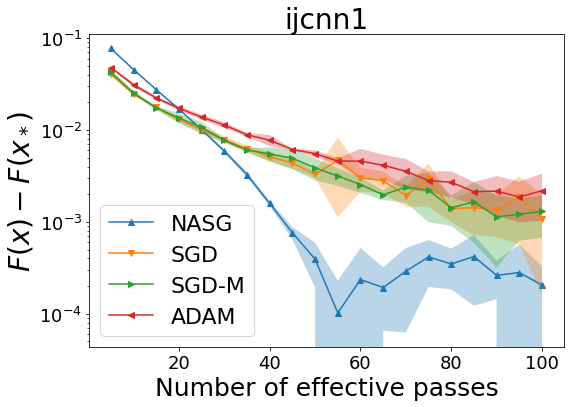}
\includegraphics[width=0.33\textwidth]{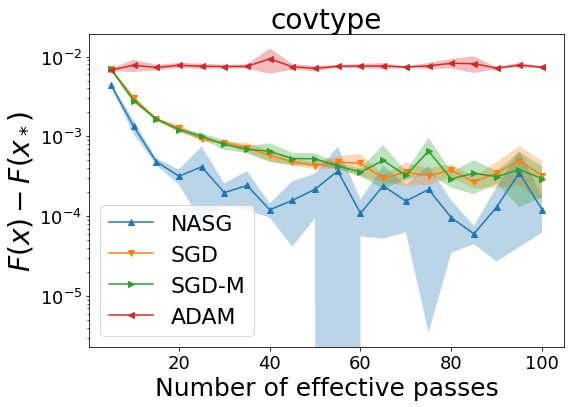}
\includegraphics[width=0.33\textwidth]{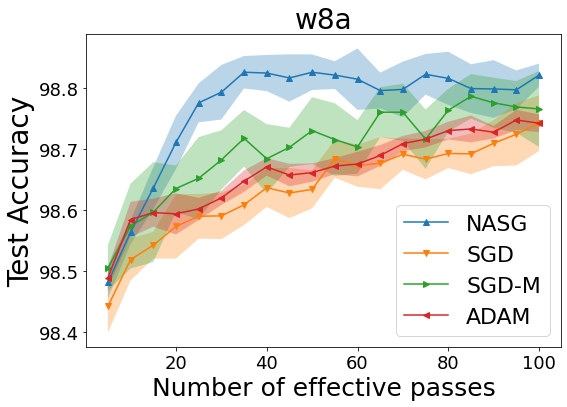}
\includegraphics[width=0.33\textwidth]{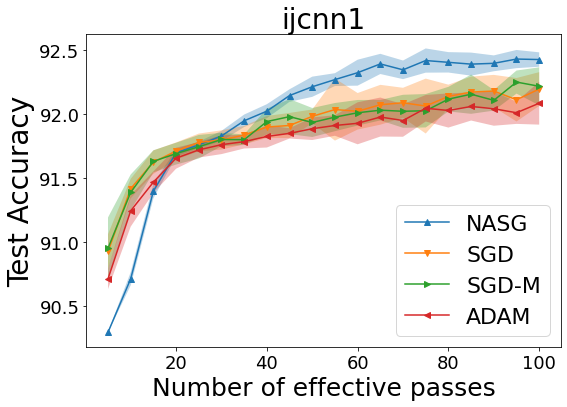}
\includegraphics[width=0.33\textwidth]{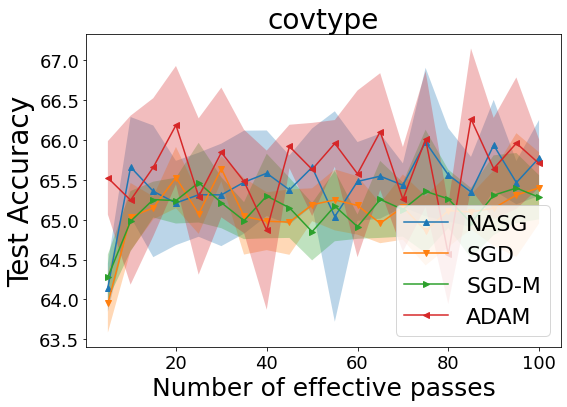}
\caption{\textbf{(Convex binary setting).} Comparisons of loss residual $F(x) - F(x_*)$ (top) and test accuracy (bottom) produced by first-order methods for \texttt{w8a}, \texttt{ijcnn1} and \texttt{covtype} datasets, respectively. The number of effective passes is the number of epochs (i.e. number of data passes) in the progress.}
\label{fig_LR}
\end{figure*}

\begin{remark}\label{rem_RR}
Let us assume that $\| \tilde{x}_0 -  x_{*} \| \leq \frac{E}{\sqrt{n}}$ where $\tilde{x}_0$ be the initial point and $E>0$ be a constant. For the same conditions as in Theorem~\ref{thm_convex_RR}, i.e. component convexity with a \textbf{\textit{randomized scheme}},
we have
\begin{align}
    \E \left[ F ( \tilde{x}_{T} ) - F ( x_{*} ) \right]
    \leq \Ocal\left(\frac{ \sigma_*^2/L + LE^2  }{nT} \right),
    \label{eq_rem_RR}
\end{align}
which shows a further improvement of $n$ over the standard setting thanks to the application of randomized schemes and initial assumption.
\end{remark}

\section{Numerical Experiments}\label{sec:experiments}

To support our theoretical analysis, we present three sets of numerical experiments, comparing our algorithm with the state-of-the-art SGD-type and shuffling gradient methods.



\subsection{Binary Classification}\label{subsec:exp_binary}
In this section, we describe the setting of Figure \ref{fig_demo} and other  experiments.
Let us consider the following convex binary classification problem:
\begin{align*}
   \min_{w \in \mathbb{R}^d} \Big \{ F(w) \! := \! \frac{1}{n} \sum_{i=1}^n \log(1 \! + \! \exp(- y_i x_i^\top w )) \!  \Big \}, 
\end{align*}
where $\sets{(x_i, y_i)}_{i=1}^n$ is a set of training samples with $x_i \in \mathbb{R}^d$ and $y_i \in \{-1,1\}$. 
We have conducted the experiments on three classification datasets \texttt{w8a} ($49,749$ samples), \texttt{ijcnn1} ($91,701$ samples) and \texttt{covtype} ($406709$ samples) from \texttt{LIBSVM} \citep{LIBSVM}. 
The stochastic experiments are repeated with random seeds 10 times and we report the average results with confidence intervals. 

In Figure \ref{fig_demo}, we compare our algorithm with NAG and the stochastic shuffling version that update Nesterov's momentum per iteration. Note that NAG is a deterministic algorithm, hence it does not have confidence intervals. In order to make fair comparisons, we report the results of three methods in Figure \ref{fig_demo} after every effective data passes, (i.e. comparing them with the same computational cost). In addition, since NAG converges slowly when $n$ is large, in our futher experiments, we choose to compare our algorithm with other stochastic first-order methods. 

In Figure \ref{fig_LR}, we compare our method with Stochastic Gradient Descent (SGD) and two stochastic algorithms: SGD with Momentum (SGD-M) \citep{Polyak1964} and Adam \citep{Kingma2014}. 
For the latter two algorithms, we use the hyper-parameter settings recommended and widely used in practice (i.e. momentum: 0.9 for SGD-M, and two hyper-parameters $\beta_1 := 0.9$, $\beta_2 := 0.999$ for Adam). 

To have a fair comparison, we apply the randomized reshuffling  scheme to all methods. 
Note that shuffling strategies are favorable in practice and have been implemented in TensorFlow, PyTorch, and Keras \citep{tensorflow2015-whitepaper, pytorch, chollet2015keras}. 
We tune each algorithm using constant learning rate and report the best final results. 

For \texttt{w8a} and \texttt{covtype} datasets, our algorithm shows better performance than the other methods in the training process. For \texttt{ijcnn1}, NASG is somewhat worse than the other methods at the beginning, however, it surpasses all other methods after a few epochs and maintains a better decrease toward the end of training stage. 
In terms of test accuracy, our method shows comparable performance for \texttt{covtype} dataset, and achieves good generalization for \texttt{w8a} and \texttt{ijcnn1} datasets.

\begin{figure*}[ht!] 
\centering
\includegraphics[width=0.33\textwidth]{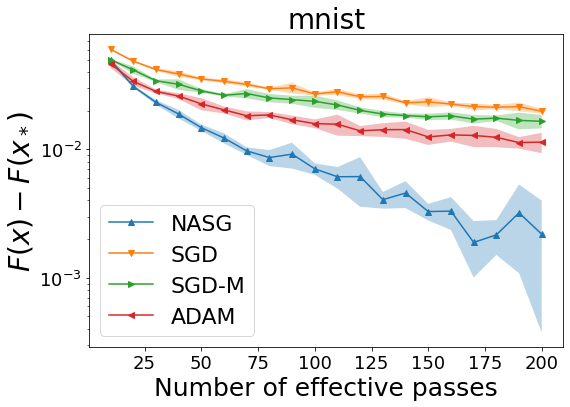}
\includegraphics[width=0.33\textwidth]{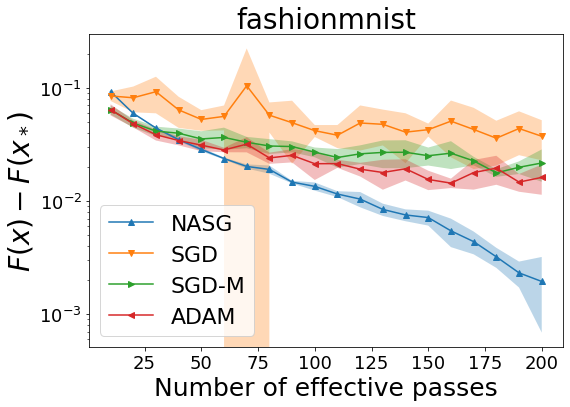}
\includegraphics[width=0.33\textwidth]{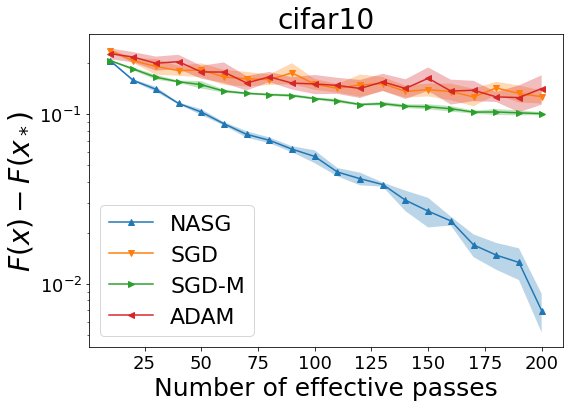}
\includegraphics[width=0.33\textwidth]{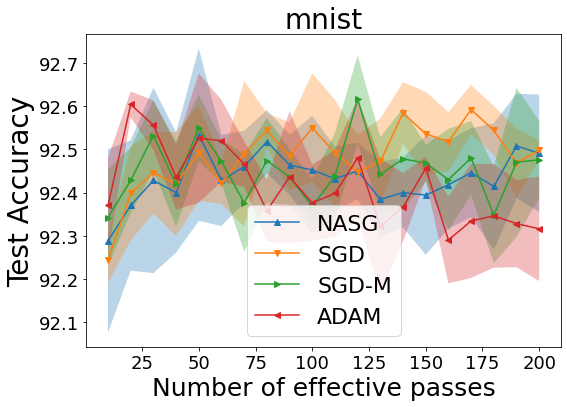}
\includegraphics[width=0.33\textwidth]{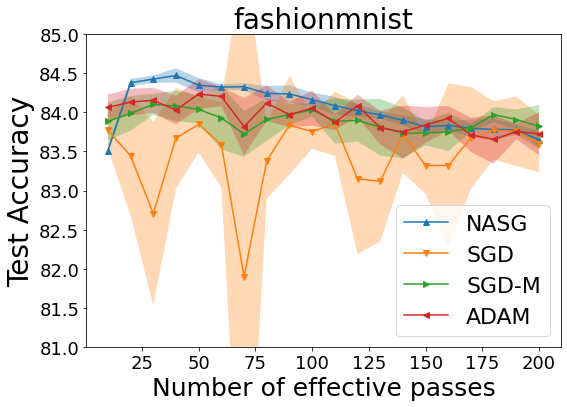}
\includegraphics[width=0.33\textwidth]{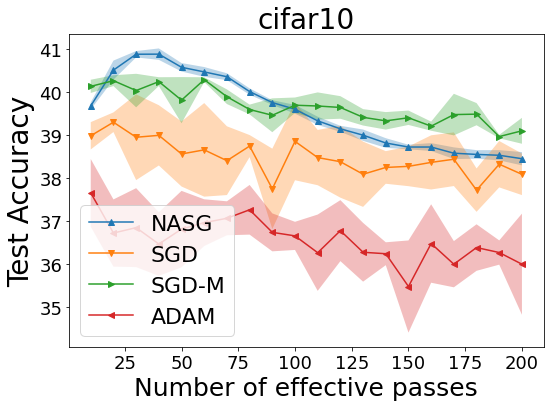}
\caption{\textbf{(Convex image setting).} Comparisons of loss residual $F(x) - F(x_*)$ (top) and test accuracy (bottom) produced by first-order methods for \texttt{MNIST}, \texttt{Fashion-MNIST} and \texttt{CIFAR-10}, respectively. The number of effective passes is the number of epochs (i.e. number of data passes) throughout the progress.}
\label{fig_NN}
\end{figure*}
In the next subsection, we perform another set of experiments in convex setting. 
Our Appendix describes all the experimental details and implementation. 



\subsection{Convex Image Classification} \label{subsec_exp_nn}

For the second experiment, we test our algorithm using linear neural networks on three well-known image classification datasets: \texttt{MNIST} dataset \citep{MNIST} and  \texttt{Fashion-MNIST} dataset \citep{xiao2017fashion} both with $60,000$ samples, and finally  \texttt{CIFAR-10} dataset \citep{CIFAR10} with $50,000$ images. 
We experiment with the following minimization problem:
\begin{align*}
    \min_{w \in \mathbb{R}^d} \Big \{F(w):= - \frac{1}{n}\sum_{i=1}^n y_i^{\top} \log( \text{softmax} ( h(w ; i) ) ) \Big \},
\end{align*}
where $h(w;i) = W x_i + b$ is a simple neural network with parameter $w = \{W,b\},$ $W \in \mathbb{R}^{c \times d}$ and $ b \in \mathbb{R}^c$. 
The input data $\sets{x_i}_{i=1}^n$ are in $\R^d$ and the output labels $\sets{y_i}_{i=1}^n$ are one-hot vectors in $\R^c$, where $c$ is the number of classes. 
The \textit{softmax} function is defined as $$\text{softmax}(z) = \left(\frac{e^{z_1}}{ \sum_{k=1}^c e^{z_k}}; \dots ; \frac{e^{z_c}}{ \sum_{k=1}^c e^{z_k}} \right)^\top. $$ 
Similar to the previous experiment, we compare our algorithm with other stochastic first-order methods with randomized reshuffling scheme. The minibatch size is 256. All the algorithms are implemented in Python using PyTorch package \citep{pytorch}. 
They are tested using 10 different random seeds and we report the average results with confidence intervals. 
We tune each algorithm using constant learning rate and report the best final results in Figure \ref{fig_NN}. 

Our algorithm achieves a better decrease than other methods on \texttt{MNIST} and \texttt{CIFAR-10} datasets very early in the training process. On \texttt{Fashion-MNIST} dataset, NASG starts slower than other methods at the beginning. In the next stage, it suggests a better performance with a little oscillations in the end. 

In terms of generalization, our method shows comparable performance to all other stochastic algorithms. Note that our main focus is the training task, that is, solving the optimization problem \eqref{ERM_problem_01} and there may be over-fitting that leads to test accuracy decrease in the later part of the training progress. 

\subsection{Non-convex Image Classification} 
We further test our algorithm with a simple non-convex model to demonstrate the efficiency and flexibility of our method beyond the convex setting. 
For this experiment, we use a similar problem as in the previous section (i.e. 
training neural networks on three image classification datasets: \texttt{MNIST},  \texttt{Fashion-MNIST} dataset  and \texttt{CIFAR-10} dataset). 
The minimization problem is:
\begin{align*}
    \min_{w \in \mathbb{R}^d} \Big \{F(w):= - \frac{1}{n}\sum_{i=1}^n y_i^{\top} \log( \text{softmax} ( h(w ; i) ) ) \Big \},
\end{align*}
where $h(w;i) = W_2 (W_1 x_i + b_1) + b_2$ is a neural network with one hidden layer containing $m$ neurons and no activation. 
The input data $\sets{x_i}_{i=1}^n$ are in $\R^d$ and the output labels $\sets{y_i}_{i=1}^n$ are one-hot vectors in $\R^c$, where $c$ is the number of classes. 
The parameter is $w = \{W_1,b_1, W_2,b_2\}$ with $W_1 \in \mathbb{R}^{m \times d}, b_1 \in \mathbb{R}^m$ and $W_2 \in \mathbb{R}^{c \times m}, b_2 \in \mathbb{R}^c$.

\begin{figure*}[ht] 
\centering
\includegraphics[width=0.33\textwidth]{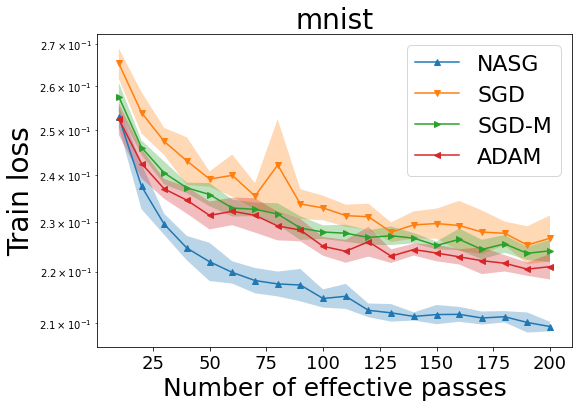}
\includegraphics[width=0.33\textwidth]{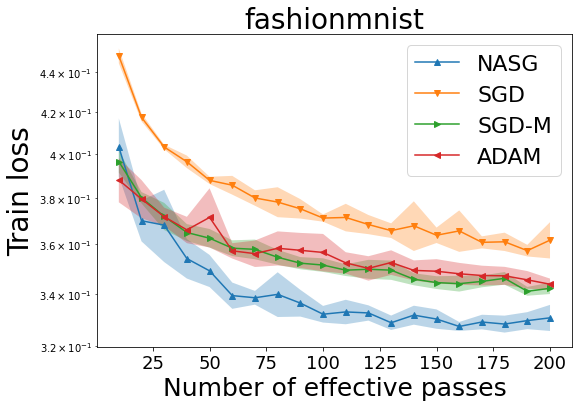}
\includegraphics[width=0.33\textwidth]{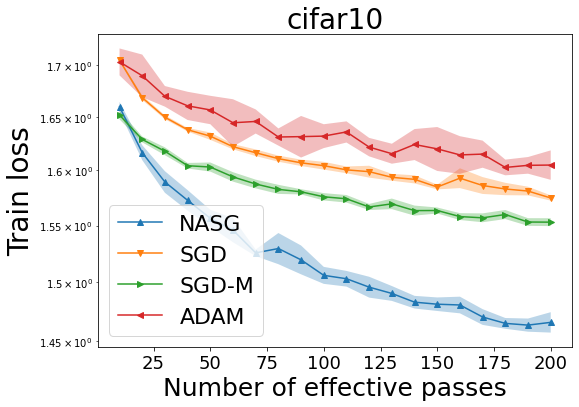}
\includegraphics[width=0.33\textwidth]{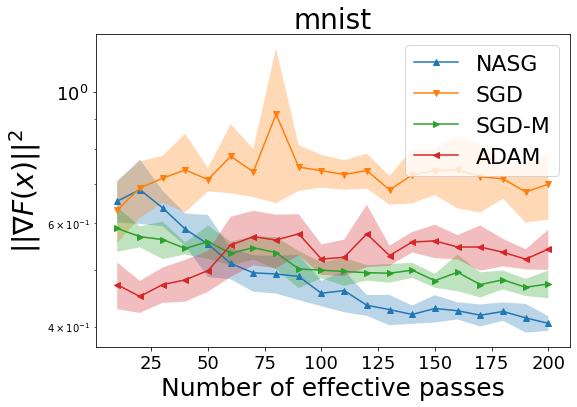}
\includegraphics[width=0.33\textwidth]{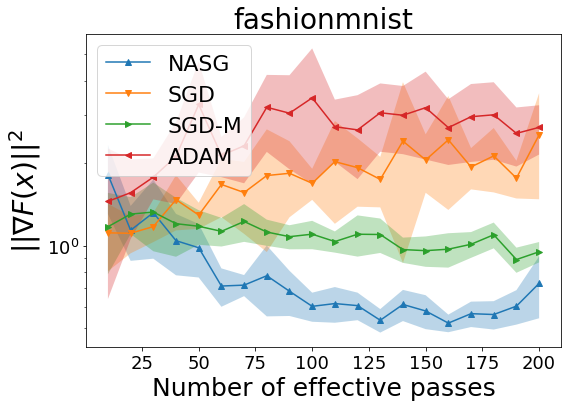}
\includegraphics[width=0.33\textwidth]{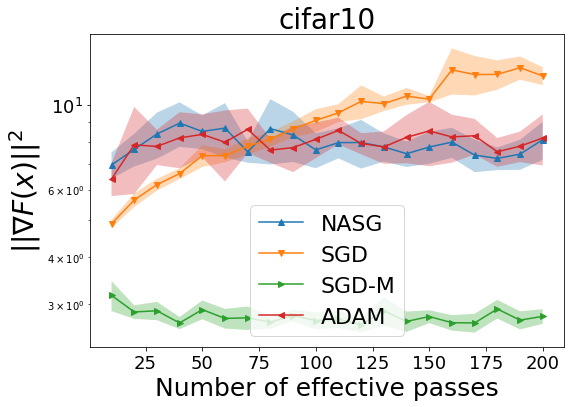}
\caption{\textbf{(Non-convex image setting).} Comparisons of train loss $F(x)$ (top) and the squared norm of gradient $\|\nabla F(x)\|^2$  produced by first-order methods for \texttt{MNIST}, \texttt{Fashion-MNIST} and \texttt{CIFAR-10}, respectively. The number of effective passes is the number of epochs (i.e number of data passes) throughout the progress.}
\label{fig_NN2}
\end{figure*}

For \texttt{MNIST} and \texttt{Fashion-MNIST} datasets, we run a small network with $m=300$ hidden neurons. For \texttt{CIFAR-10} dataset, we experiment with $m=900$ neurons in the hidden layer.
Similar to the previous experiments, we compare our algorithm with other stochastic first-order methods (SGD, SGD-M and ADAM) and we apply randomized reshuffling scheme to all these algorithms. We implement all the methods in Python using PyTorch package \citep{pytorch}, then tune each algorithm using constant learning rate. Figure \ref{fig_NN2} report the train loss and the squared norm of gradient returned by our experiments. We delay other experimental setting details to the Appendix\footnote{Our code can be found at the repository \url{https://github.com/htt-trangtran/nasg}.}.

In this simple non-convex setting, our algorithm also achieves a better decrease than other methods on \texttt{MNIST} and \texttt{CIFAR-10} datasets. On \texttt{Fashion-MNIST} dataset, NASG starts slower than other methods at the beginning and surpasses other method later in the training process. 
In terms of gradient norm, our method shows competitive performance for \texttt{MNIST} and \texttt{Fashion-MNIST} datasets, while performs comparably good in \texttt{CIFAR-10} dataset. We further note that all our experiments are tuned to the best value of the training loss for every algorithm. 


\section{Conclusions and Future Work}\label{sec_conclusion}

We propose Nesterov Accelerated Shuffling Gradient (NASG), a new gradient method that combines the update of SGD using shuffling sampling schemes with Nesterov's momentum. 
Our method achieves a convergence rate of $\Ocal(1/T)$ for smooth convex functions, where $T$ is the number of effective data passes. This rate is better than the state-of-the-art result of SGD using shuffling schemes, in terms of $T$. 

Although we have made progresses in understanding theoretical properties of shuffling methods in general (and NASG in particular), an interesting research question remains: whether our method can achieve a better theoretical rate in terms of the number of data points $n$. Our work answers this question partially by different approaches, including the application of randomized sampling schemes and the investigation of an initial condition. 
In addition, investigating our algorithm in non-convex settings is also a promising direction. 

\section*{Acknowledgements}

The authors would like to thank the reviewers for their useful
comments and suggestions which helped to improve the
exposition of this paper. The work of Trang H. Tran and Katya Scheinberg have partly been supported by the ONR Grant N00014-22-1-2154 and the NSF Grant CCF 21-40057.

\newpage
\bibliography{reference}

\begin{thebibliography}{52}
\providecommand{\natexlab}[1]{#1}
\providecommand{\url}[1]{\texttt{#1}}
\expandafter\ifx\csname urlstyle\endcsname\relax
  \providecommand{\doi}[1]{doi: #1}\else
  \providecommand{\doi}{doi: \begingroup \urlstyle{rm}\Url}\fi

\bibitem[Abadi et~al.(2015)Abadi, Agarwal, Barham, Brevdo, Chen, Citro,
  Corrado, Davis, Dean, Devin, Ghemawat, Goodfellow, Harp, Irving, Isard, Jia,
  Jozefowicz, Kaiser, Kudlur, Levenberg, Man\'{e}, Monga, Moore, Murray, Olah,
  Schuster, Shlens, Steiner, Sutskever, Talwar, Tucker, Vanhoucke, Vasudevan,
  Vi\'{e}gas, Vinyals, Warden, Wattenberg, Wicke, Yu, and
  Zheng]{tensorflow2015-whitepaper}
Abadi, M., Agarwal, A., Barham, P., Brevdo, E., Chen, Z., Citro, C., Corrado,
  G.~S., Davis, A., Dean, J., Devin, M., Ghemawat, S., Goodfellow, I., Harp,
  A., Irving, G., Isard, M., Jia, Y., Jozefowicz, R., Kaiser, L., Kudlur, M.,
  Levenberg, J., Man\'{e}, D., Monga, R., Moore, S., Murray, D., Olah, C.,
  Schuster, M., Shlens, J., Steiner, B., Sutskever, I., Talwar, K., Tucker, P.,
  Vanhoucke, V., Vasudevan, V., Vi\'{e}gas, F., Vinyals, O., Warden, P.,
  Wattenberg, M., Wicke, M., Yu, Y., and Zheng, X.
\newblock {TensorFlow}: Large-scale machine learning on heterogeneous systems,
  2015.
\newblock URL \url{https://www.tensorflow.org/}.
\newblock Software available from tensorflow.org.

\bibitem[Ahn et~al.(2020)Ahn, Yun, and Sra]{ahn2020sgd}
Ahn, K., Yun, C., and Sra, S.
\newblock Sgd with shuffling: optimal rates without component convexity and
  large epoch requirements.
\newblock \emph{arXiv preprint arXiv:2006.06946}, 2020.

\bibitem[Bottou(2009)]{bottou2009curiously}
Bottou, L.
\newblock Curiously fast convergence of some stochastic gradient descent
  algorithms.
\newblock In \emph{Proceedings of the symposium on learning and data science,
  Paris}, volume~8, pp.\  2624--2633, 2009.

\bibitem[Bottou(2012)]{bottou2012stochastic}
Bottou, L.
\newblock Stochastic gradient descent tricks.
\newblock In \emph{Neural networks: Tricks of the trade}, pp.\  421--436.
  Springer, 2012.

\bibitem[Bottou et~al.(2018)Bottou, Curtis, and Nocedal]{Bottou2018}
Bottou, L., Curtis, F.~E., and Nocedal, J.
\newblock {O}ptimization {M}ethods for {L}arge-{S}cale {M}achine {L}earning.
\newblock \emph{SIAM Rev.}, 60\penalty0 (2):\penalty0 223--311, 2018.

\bibitem[Chang \& Lin(2011)Chang and Lin]{LIBSVM}
Chang, C.-C. and Lin, C.-J.
\newblock {LIBSVM}: A library for support vector machines.
\newblock \emph{ACM Transactions on Intelligent Systems and Technology},
  2:\penalty0 27:1--27:27, 2011.
\newblock Software available at \url{http://www.csie.ntu.edu.tw/~cjlin/libsvm}.

\bibitem[Chollet et~al.(2015)]{chollet2015keras}
Chollet, F. et~al.
\newblock Keras.
\newblock \emph{GitHub}, 2015.
\newblock URL \url{https://github.com/fchollet/keras}.

\bibitem[Defazio et~al.(2014)Defazio, Bach, and Lacoste-Julien]{SAGA}
Defazio, A., Bach, F., and Lacoste-Julien, S.
\newblock Saga: A fast incremental gradient method with support for
  non-strongly convex composite objectives.
\newblock In \emph{Advances in Neural Information Processing Systems}, pp.\
  1646--1654, 2014.

\bibitem[Devolder et~al.(2014)Devolder, Glineur, and
  Nesterov]{springer_devolder2014}
Devolder, O., Glineur, F., and Nesterov, Y.
\newblock First-order methods of smooth convex optimization with inexact
  oracle.
\newblock \emph{Springer-Verlag}, 146\penalty0 (1–2), 2014.
\newblock ISSN 0025-5610.
\newblock \doi{10.1007/s10107-013-0677-5}.
\newblock URL \url{https://doi.org/10.1007/s10107-013-0677-5}.

\bibitem[Dozat(2016)]{dozat2016incorporating}
Dozat, T.
\newblock Incorporating nesterov momentum into {ADAM}.
\newblock \emph{ICLR Workshop}, 1:\penalty0 2013–--2016, 2016.

\bibitem[Duchi et~al.(2011)Duchi, Hazan, and Singer]{AdaGrad}
Duchi, J., Hazan, E., and Singer, Y.
\newblock Adaptive subgradient methods for online learning and stochastic
  optimization.
\newblock \emph{Journal of Machine Learning Research}, 12:\penalty0 2121--2159,
  2011.

\bibitem[Ghadimi \& Lan(2013)Ghadimi and Lan]{ghadimi2013stochastic}
Ghadimi, S. and Lan, G.
\newblock Stochastic first-and zeroth-order methods for nonconvex stochastic
  programming.
\newblock \emph{SIAM J. Optim.}, 23\penalty0 (4):\penalty0 2341--2368, 2013.

\bibitem[Gürbüzbalaban et~al.(2019)Gürbüzbalaban, Ozdaglar, and
  Parrilo]{Gurbuzbalaban2019}
Gürbüzbalaban, M., Ozdaglar, A., and Parrilo, P.~A.
\newblock Why random reshuffling beats stochastic gradient descent.
\newblock \emph{Mathematical Programming}, Oct 2019.
\newblock ISSN 1436-4646.
\newblock \doi{10.1007/s10107-019-01440-w}.

\bibitem[Haochen \& Sra(2019)Haochen and Sra]{haochen2019random}
Haochen, J. and Sra, S.
\newblock Random shuffling beats sgd after finite epochs.
\newblock In \emph{International Conference on Machine Learning}, pp.\
  2624--2633. PMLR, 2019.

\bibitem[Hu et~al.(2009)Hu, Pan, and Kwok]{NIPS_hu2009}
Hu, C., Pan, W., and Kwok, J.
\newblock Accelerated gradient methods for stochastic optimization and online
  learning.
\newblock In Bengio, Y., Schuurmans, D., Lafferty, J., Williams, C., and
  Culotta, A. (eds.), \emph{Advances in Neural Information Processing Systems},
  volume~22. Curran Associates, Inc., 2009.
\newblock URL
  \url{https://proceedings.neurips.cc/paper/2009/file/ec5aa0b7846082a2415f0902f0da88f2-Paper.pdf}.

\bibitem[Johnson \& Zhang(2013)Johnson and Zhang]{SVRG}
Johnson, R. and Zhang, T.
\newblock Accelerating stochastic gradient descent using predictive variance
  reduction.
\newblock In \emph{NIPS}, pp.\  315--323, 2013.

\bibitem[Kingma \& Ba(2014)Kingma and Ba]{Kingma2014}
Kingma, D.~P. and Ba, J.
\newblock {ADAM}: {A} {M}ethod for {S}tochastic {O}ptimization.
\newblock \emph{Proceedings of the 3rd International Conference on Learning
  Representations (ICLR)}, abs/1412.6980, 2014.

\bibitem[Krizhevsky \& Hinton(2009)Krizhevsky and Hinton]{CIFAR10}
Krizhevsky, A. and Hinton, G.
\newblock Learning multiple layers of features from tiny images.
\newblock Technical report, Citeseer, 2009.

\bibitem[Lan(2012)]{Lan2012}
Lan, G.
\newblock An optimal method for stochastic composite optimization.
\newblock \emph{Mathematical Programming}, 133:\penalty0 365--397, 2012.

\bibitem[Le~Roux et~al.(2012)Le~Roux, Schmidt, and Bach]{SAG}
Le~Roux, N., Schmidt, M., and Bach, F.
\newblock A stochastic gradient method with an exponential convergence rate for
  finite training sets.
\newblock In \emph{NIPS}, pp.\  2663--2671, 2012.

\bibitem[LeCun et~al.(1998)LeCun, Bottou, Bengio, and Haffner]{MNIST}
LeCun, Y., Bottou, L., Bengio, Y., and Haffner, P.
\newblock Gradient-based learning applied to document recognition.
\newblock \emph{Proceedings of the IEEE}, 86\penalty0 (11):\penalty0
  2278--2324, 1998.

\bibitem[Lessard et~al.(2016)Lessard, Recht, and Packard]{siam_lessard2016}
Lessard, L., Recht, B., and Packard, A.
\newblock Analysis and design of optimization algorithms via integral quadratic
  constraints.
\newblock \emph{SIAM Journal on Optimization}, 26\penalty0 (1):\penalty0
  57--95, 2016.
\newblock \doi{10.1137/15M1009597}.
\newblock URL \url{https://doi.org/10.1137/15M1009597}.

\bibitem[Liu \& Belkin(2018)Liu and Belkin]{corr_liu2018}
Liu, C. and Belkin, M.
\newblock Mass: an accelerated stochastic method for over-parametrized
  learning.
\newblock \emph{CoRR}, abs/1810.13395, 2018.
\newblock URL \url{http://arxiv.org/abs/1810.13395}.

\bibitem[Mishchenko et~al.(2020)Mishchenko, Khaled Ragab~Bayoumi, and
  Richt{\'a}rik]{mishchenko2020random}
Mishchenko, K., Khaled Ragab~Bayoumi, A., and Richt{\'a}rik, P.
\newblock Random reshuffling: Simple analysis with vast improvements.
\newblock \emph{Advances in Neural Information Processing Systems}, 33, 2020.

\bibitem[Mishchenko et~al.(2021)Mishchenko, Khaled, and
  Richtárik]{mishchenko2021proximal}
Mishchenko, K., Khaled, A., and Richtárik, P.
\newblock Proximal and federated random reshuffling, 2021.

\bibitem[Nagaraj et~al.(2019)Nagaraj, Jain, and Netrapalli]{nagaraj2019sgd}
Nagaraj, D., Jain, P., and Netrapalli, P.
\newblock Sgd without replacement: Sharper rates for general smooth convex
  functions.
\newblock In \emph{International Conference on Machine Learning}, pp.\
  4703--4711, 2019.

\bibitem[Nedic \& Bertsekas(2001{\natexlab{a}})Nedic and
  Bertsekas]{nedic2001convergence}
Nedic, A. and Bertsekas, D.
\newblock Convergence rate of incremental subgradient algorithms.
\newblock In \emph{Stochastic optimization: algorithms and applications}, pp.\
  223--264. Springer, 2001{\natexlab{a}}.

\bibitem[Nedic \& Bertsekas(2001{\natexlab{b}})Nedic and
  Bertsekas]{nedic2001incremental}
Nedic, A. and Bertsekas, D.~P.
\newblock Incremental subgradient methods for nondifferentiable optimization.
\newblock \emph{SIAM J. on Optim.}, 12\penalty0 (1):\penalty0 109--138,
  2001{\natexlab{b}}.

\bibitem[Nemirovski et~al.(2009)Nemirovski, Juditsky, Lan, and
  Shapiro]{Nemirovski2009}
Nemirovski, A., Juditsky, A., Lan, G., and Shapiro, A.
\newblock Robust stochastic approximation approach to stochastic programming.
\newblock \emph{SIAM J. on Optimization}, 19\penalty0 (4):\penalty0 1574--1609,
  2009.

\bibitem[Nesterov(1983)]{Nesterov1983}
Nesterov, Y.
\newblock A method for unconstrained convex minimization problem with the rate
  of convergence $\mathcal{O}(1/k^2)$.
\newblock \emph{Doklady AN SSSR}, 269:\penalty0 543--547, 1983.
\newblock Translated as Soviet Math. Dokl.

\bibitem[Nesterov(2004)]{Nesterov2004}
Nesterov, Y.
\newblock \emph{{I}ntroductory lectures on convex optimization: {A} basic
  course}, volume~87 of \emph{Applied Optimization}.
\newblock Kluwer Academic Publishers, 2004.

\bibitem[Nguyen et~al.(2018)Nguyen, Nguyen, van Dijk, Richtarik, Scheinberg,
  and Takac]{Nguyen2018_sgdhogwild}
Nguyen, L., Nguyen, P.~H., van Dijk, M., Richtarik, P., Scheinberg, K., and
  Takac, M.
\newblock {SGD} and {H}ogwild! convergence without the bounded gradients
  assumption.
\newblock In \emph{Proceedings of the 35th International Conference on Machine
  Learning-Volume 80}, pp.\  3747--3755, 2018.

\bibitem[Nguyen et~al.(2017)Nguyen, Liu, Scheinberg, and
  Tak{\'a}{\v{c}}]{Nguyen2017sarah}
Nguyen, L.~M., Liu, J., Scheinberg, K., and Tak{\'a}{\v{c}}, M.
\newblock {SARAH}: A novel method for machine learning problems using
  stochastic recursive gradient.
\newblock In \emph{Proceedings of the 34th International Conference on Machine
  Learning-Volume 70}, pp.\  2613--2621. JMLR. org, 2017.

\bibitem[Nguyen et~al.(2021)Nguyen, Tran-Dinh, Phan, Nguyen, and van
  Dijk]{nguyen2020unified}
Nguyen, L.~M., Tran-Dinh, Q., Phan, D.~T., Nguyen, P.~H., and van Dijk, M.
\newblock A unified convergence analysis for shuffling-type gradient methods.
\newblock \emph{Journal of Machine Learning Research}, 22\penalty0
  (207):\penalty0 1--44, 2021.

\bibitem[Nguyen et~al.(2022)Nguyen, Tran, and van Dijk]{Nguyen2022_OptDL}
Nguyen, L.~M., Tran, T.~H., and van Dijk, M.
\newblock Finite-sum optimization: A new perspective for convergence to a
  global solution.
\newblock \emph{arXiv preprint arXiv:2202.03524}, 2022.

\bibitem[Paszke et~al.(2019)Paszke, Gross, Massa, Lerer, Bradbury, Chanan,
  Killeen, Lin, Gimelshein, Antiga, Desmaison, Kopf, Yang, DeVito, Raison,
  Tejani, Chilamkurthy, Steiner, Fang, Bai, and Chintala]{pytorch}
Paszke, A., Gross, S., Massa, F., Lerer, A., Bradbury, J., Chanan, G., Killeen,
  T., Lin, Z., Gimelshein, N., Antiga, L., Desmaison, A., Kopf, A., Yang, E.,
  DeVito, Z., Raison, M., Tejani, A., Chilamkurthy, S., Steiner, B., Fang, L.,
  Bai, J., and Chintala, S.
\newblock Pytorch: An imperative style, high-performance deep learning library.
\newblock In \emph{Advances in Neural Information Processing Systems 32}, pp.\
  8024--8035. Curran Associates, Inc., 2019.

\bibitem[Polyak \& Juditsky(1992)Polyak and Juditsky]{Polyak1992}
Polyak, B. and Juditsky, A.
\newblock Acceleration of stochastic approximation by averaging.
\newblock \emph{SIAM J. Control Optim.}, 30\penalty0 (4):\penalty0 838--855,
  1992.

\bibitem[Polyak(1964)]{Polyak1964}
Polyak, B.~T.
\newblock Some methods of speeding up the convergence of iteration methods.
\newblock \emph{{USSR} Computational Mathematics and Mathematical Physics},
  4\penalty0 (5):\penalty0 1--17, 1964.

\bibitem[Rajput et~al.(2020)Rajput, Gupta, and
  Papailiopoulos]{rajput2020closing}
Rajput, S., Gupta, A., and Papailiopoulos, D.
\newblock Closing the convergence gap of sgd without replacement.
\newblock In \emph{International Conference on Machine Learning}, pp.\
  7964--7973. PMLR, 2020.

\bibitem[Recht \& Ré(2011)Recht and Ré]{Recht2011}
Recht, B. and Ré, C.
\newblock Parallel stochastic gradient algorithms for large-scale matrix
  completion.
\newblock \emph{Mathematical Programming Computation}, 5, 04 2011.
\newblock \doi{10.1007/s12532-013-0053-8}.

\bibitem[Robbins \& Monro(1951)Robbins and Monro]{RM1951}
Robbins, H. and Monro, S.
\newblock A stochastic approximation method.
\newblock \emph{The Annals of Mathematical Statistics}, 22\penalty0
  (3):\penalty0 400--407, 1951.

\bibitem[Safran \& Shamir(2020)Safran and Shamir]{safran2020good}
Safran, I. and Shamir, O.
\newblock How good is sgd with random shuffling?
\newblock In \emph{Conference on Learning Theory}, pp.\  3250--3284. PMLR,
  2020.

\bibitem[Schmidt \& Roux(2013)Schmidt and Roux]{schmidt2013fast}
Schmidt, M. and Roux, N.~L.
\newblock Fast convergence of stochastic gradient descent under a strong growth
  condition, 2013.

\bibitem[Shamir(2016)]{shamir2016without}
Shamir, O.
\newblock Without-replacement sampling for stochastic gradient methods.
\newblock In \emph{Advances in neural information processing systems}, pp.\
  46--54, 2016.

\bibitem[Shamir \& Zhang(2013)Shamir and Zhang]{pmlr-v28-shamir13}
Shamir, O. and Zhang, T.
\newblock Stochastic gradient descent for non-smooth optimization: Convergence
  results and optimal averaging schemes.
\newblock In Dasgupta, S. and McAllester, D. (eds.), \emph{Proceedings of the
  30th International Conference on Machine Learning}, volume~28 of
  \emph{Proceedings of Machine Learning Research}, pp.\  71--79, Atlanta,
  Georgia, USA, 17--19 Jun 2013. PMLR.
\newblock URL \url{https://proceedings.mlr.press/v28/shamir13.html}.

\bibitem[Sra et~al.(2012)Sra, Nowozin, and Wright]{sra2012optimization}
Sra, S., Nowozin, S., and Wright, S.~J.
\newblock \emph{{O}ptimization for {M}achine {L}earning}.
\newblock MIT Press, 2012.

\bibitem[Sutskever et~al.(2013)Sutskever, Martens, Dahl, and
  Hinton]{pmlr-v28-sutskever13}
Sutskever, I., Martens, J., Dahl, G., and Hinton, G.
\newblock On the importance of initialization and momentum in deep learning.
\newblock In Dasgupta, S. and McAllester, D. (eds.), \emph{Proceedings of the
  30th International Conference on Machine Learning}, volume~28 of
  \emph{Proceedings of Machine Learning Research}, pp.\  1139--1147, Atlanta,
  Georgia, USA, 17--19 Jun 2013. PMLR.
\newblock URL \url{https://proceedings.mlr.press/v28/sutskever13.html}.

\bibitem[Tran et~al.(2021)Tran, Nguyen, and Tran-Dinh]{pmlr-v139-tran21b}
Tran, T.~H., Nguyen, L.~M., and Tran-Dinh, Q.
\newblock {SMG}: A shuffling gradient-based method with momentum.
\newblock In Meila, M. and Zhang, T. (eds.), \emph{Proceedings of the 38th
  International Conference on Machine Learning}, volume 139 of
  \emph{Proceedings of Machine Learning Research}, pp.\  10379--10389. PMLR,
  18--24 Jul 2021.
\newblock URL \url{https://proceedings.mlr.press/v139/tran21b.html}.

\bibitem[Vaswani et~al.(2019)Vaswani, Bach, and Schmidt]{pmlr-v89-vaswani19a}
Vaswani, S., Bach, F., and Schmidt, M.
\newblock Fast and faster convergence of sgd for over-parameterized models and
  an accelerated perceptron.
\newblock In Chaudhuri, K. and Sugiyama, M. (eds.), \emph{Proceedings of the
  Twenty-Second International Conference on Artificial Intelligence and
  Statistics}, volume~89 of \emph{Proceedings of Machine Learning Research},
  pp.\  1195--1204. PMLR, 16--18 Apr 2019.
\newblock URL \url{https://proceedings.mlr.press/v89/vaswani19a.html}.

\bibitem[Xiao et~al.(2017)Xiao, Rasul, and Vollgraf]{xiao2017fashion}
Xiao, H., Rasul, K., and Vollgraf, R.
\newblock Fashion-mnist: a novel image dataset for benchmarking machine
  learning algorithms.
\newblock \emph{arXiv preprint arXiv:1708.07747}, 2017.

\bibitem[Yuan et~al.(2016)Yuan, Ying, and Sayed]{JMLR_yuan2016}
Yuan, K., Ying, B., and Sayed, A.~H.
\newblock On the influence of momentum acceleration on online learning.
\newblock \emph{Journal of Machine Learning Research}, 17\penalty0
  (192):\penalty0 1--66, 2016.
\newblock URL \url{http://jmlr.org/papers/v17/16-157.html}.

\bibitem[Zhong \& Kwok(2014)Zhong and Kwok]{pmlr-v33-zhong14}
Zhong, W. and Kwok, J.
\newblock {Accelerated Stochastic Gradient Method for Composite
  Regularization}.
\newblock In Kaski, S. and Corander, J. (eds.), \emph{Proceedings of the
  Seventeenth International Conference on Artificial Intelligence and
  Statistics}, volume~33 of \emph{Proceedings of Machine Learning Research},
  pp.\  1086--1094, Reykjavik, Iceland, 22--25 Apr 2014. PMLR.
\newblock URL \url{https://proceedings.mlr.press/v33/zhong14.html}.

\end{thebibliography}
\bibliographystyle{icml2022}

\newpage
\appendix
\onecolumn


  \vbox{%
    \hsize\textwidth
    \linewidth\hsize
    \vskip 0.1in
  \hrule height 4pt
  \vskip 0.25in
  \vskip -5.5pt%
  \centering
    {\Large \bf{Nesterov Accelerated Shuffling Gradient Method for Convex Optimization \\
    Appendix, ICML 2022} \par}
      \vskip 0.29in
  \vskip -5.5pt
  \hrule height 1pt
  \vskip 0.09in%
    
  \vskip 0.2in
  }

\appendix
\section{Technical Lemmas}\label{sec_app_lemma}
\subsection{Basic Derivations for Algorithm \ref{shuffling_nesterov_02}}
Let us collect all the basic necessary expressions for Algorithm \ref{shuffling_nesterov_02}. From the update $y_{i}^{(t)} := y_{i-1}^{(t)} - \eta_i^{(t)} \nabla f ( y_{i-1}^{(t)} ; \pi^{(t)} ( i ) )$, we have the following for $i = 1,\dots,n, t \geq 1$:
\begin{align*}
    y_{i}^{(t)} = y_{i-1}^{(t)} - \eta_i^{(t)} \nabla f ( y_{i-1}^{(t)} ; \pi^{(t)} ( i ) ) = y_{0}^{(t)} - \sum_{j=1}^{i} \eta_j^{(t)} \nabla f ( y_{j-1}^{(t)} ; \pi^{(t)} ( j ) ). \tagthis \label{eq_003}
\end{align*}

Note that $y_0^{(t)} = \tilde{y}_{t-1}$ and $\tilde{x}_t = y_{n}^{(t)}$ and $\eta_i^{(t)} = \frac{\eta_t}{n}$ for $i = 1,\dots,n, t \geq 1$, we have
\begin{align*}
    \tilde{x}_t = \tilde{y}_{t-1} - \frac{\eta_t}{n} \sum_{j=1}^{n} \nabla f ( y_{j-1}^{(t)} ; \pi^{(t)} ( j ) ). \tagthis \label{eq_004}
\end{align*}



From \eqref{eq_006_0} and \eqref{eq_006} we have the following for $t \geq 1$:
\begin{align*}
    \tilde{x}_{t+1} = \tilde{y}_{t} + \theta^{(t)} (v^{(t+1)} - v^{(t)}). \tagthis \label{eq_006_2}
\end{align*}

On the other hand, we consider the term $v^{(t+1)} = \frac{t+2}{2} \tilde{x}_{t+1} - \frac{t}{2} \tilde{x}_{t}$ for $t \geq 0$:
\begin{align*}
    \frac{t+2}{2} \tilde{x}_{t+1} - \frac{t}{2} \tilde{x}_{t} & \overset{\eqref{eq_004}}{=} \frac{t+2}{2} \left( \tilde{y}_{t} - \frac{\eta_{t+1}}{n} \sum_{j=1}^{n} \nabla f ( y_{j-1}^{(t+1)} ; \pi^{(t+1)} ( j ) )  \right) - \frac{t}{2} \tilde{x}_{t} \\
    & \overset{\eqref{eq_005}}{=} \frac{t+2}{2} \left( \tilde{x}_{t} + \frac{t-1}{t+2} ( \tilde{x}_{t} - \tilde{x}_{t-1} ) - \frac{\eta_{t+1}}{n} \sum_{j=1}^{n} \nabla f ( y_{j-1}^{(t+1)} ; \pi^{(t+1)} ( j ) )  \right) - \frac{t}{2} \tilde{x}_{t} \\
    & = \left( \frac{t+1}{2} \tilde{x}_{t} - \frac{t-1}{2} \tilde{x}_{t-1} \right) - \left( \frac{t+2}{2} \right) \frac{\eta_{t+1}}{n} \sum_{j=1}^{n} \nabla f ( y_{j-1}^{(t+1)} ; \pi^{(t+1)} ( j ) ). 
\end{align*}

Therefore by definitions of $v^{(t)}$ and $\theta^{(t)}$ we have
\begin{align*}
    v^{(t+1)} &= v^{(t)} - \frac{\eta_{t+1}}{\theta^{(t)}} \cdot \frac{1}{n} \sum_{j=1}^{n} \nabla f ( y_{j-1}^{(t+1)} ; \pi^{(t+1)} ( j ) ),\  t\geq 0. \tagthis \label{eq_007}
\end{align*}

Using convexity of $F$, for any $x \in \mathbb{R}^d$, $y \in \mathbb{R}^d$, and $\theta \in [0,1]$, we have
\begin{align*}
    (1 - \theta) F ( x ) + \theta F ( x_{*} ) \geq F ( (1 - \theta) x + \theta x_{*} ) \geq F(y) + \langle \nabla F ( y ) ,  (1 - \theta) x + \theta x_{*} - y \rangle, 
\end{align*}
where $x_* = \arg \min_{x} F(x)$ is an optimal solution of $F$. Hence, 
\begin{align*}
    F(y) \leq (1 - \theta) F ( x ) + \theta F ( x_{*} ) + \langle \nabla F ( y ) ,  y - (1 - \theta) x - \theta x_{*} \rangle. \tagthis \label{eq_008}
\end{align*}
In addition, we define the following term: 
\begin{align*}
K_t &= \frac{1}{n} \sum_{i=1}^{n} \left \| y_{i}^{(t)} - y_{0}^{(t)}  \right \|^2\text{ and }
I_t = \frac{1}{n} \sum_{i=1}^{n} \left \| y_{n}^{(t)} - y_{i}^{(t)}  \right \|^2.
\end{align*}

For each epoch $t=1, \cdots, T$, we denote  $\mathcal{F}_t $ by $ \sigma(y_0^{(1)},\cdots,y_0^{(t)})$, the $\sigma$-algebra generated by the iterates of Algorithm~\ref{shuffling_nesterov_02} up to the beginning of the epoch $t$.
We also denote $\E_t[\cdot] $ by $\E [\cdot \mid \mathcal{F}_t]$, the conditional expectation 
on the $\sigma$-algebra $\mathcal{F}_t$.
\subsection{Key Lemmas and Proofs of Key Lemmas}

\begin{lemma}\label{lem_bound_common}
Suppose that Assumption \ref{ass_smooth} holds for \eqref{ERM_problem_01}, and $F$ is convex. 
Let $\sets{\tilde{x}_{t}}$ be generated by  Algorithm~\ref{shuffling_nesterov_02} with the learning rate $\eta_i^{(t)} := \frac{\eta_t}{n} > 0$ for a given positive sequence $\sets{\eta_t}$ with $\eta_t \leq \frac{1}{L} $. Let $\epsilon_t$ be a positive sequence, $t \geq 1$.
Then we have the following for $t \geq 1$:
\begin{align*}
    T(T+2) [F ( \tilde{x}_{T} ) - F ( x_{*} ) ]
    & \leq   \sum_{t=1}^T\frac{L^2\eta_t(t+1)^2}{2 \epsilon_t}  K_t -\sum_{t=1}^T [F ( \tilde{x}_{t} ) - F ( x_{*} ) ]\\
    & \qquad+ \sum_{t=1}^T\frac{2}{\eta_t}  \| v^{(t-1)} -  x_{*} \|^2 - \sum_{t=1}^T \frac{2}{\eta_t}(1-\epsilon_t) \| v^{(t)} -  x_{*} \|^2. \tagthis \label{eq_020_0}
\end{align*}
\end{lemma}
\subsection*{Proof of Lemma \ref{lem_bound_common}: Key estimate for Algorithm~\ref{shuffling_nesterov_02}} 
We start with the update \eqref{eq_004} of Algorithm~\ref{shuffling_nesterov_02}, for $t \geq 1$:
\allowdisplaybreaks
\begin{align*}
    F(\tilde{x}_t) & \overset{\eqref{eq_004}}{=} F \left( \tilde{y}_{t-1} - \frac{\eta_t}{n} \sum_{j=1}^{n} \nabla f ( y_{j-1}^{(t)} ; \pi^{(t)} ( j ) ) \right) \\
    & \overset{\eqref{eq:Lsmooth}}{\leq} F ( \tilde{y}_{t-1} ) - \eta_t \left \langle \nabla F ( \tilde{y}_{t-1} ) ,  \frac{1}{n} \sum_{j=1}^{n} \nabla f ( y_{j-1}^{(t)} ; \pi^{(t)} ( j ) ) \right \rangle + \frac{L \eta_t^2}{2} \left \| \frac{1}{n} \sum_{j=1}^{n} \nabla f ( y_{j-1}^{(t)} ; \pi^{(t)} ( j ) )  \right \|^2 \\
    & \overset{\eqref{eq_008},\eqref{eq_006}}{\leq} (1 - \theta^{(t-1)}) F ( \tilde{x}_{t-1} ) + \theta^{(t-1)} F ( x_{*} ) + \langle \nabla F ( \tilde{y}_{t-1} ) ,  \theta^{(t-1)} v^{(t-1)} - \theta^{(t-1)} x_{*} \rangle \\
    & \qquad -  \eta_t \left \langle \nabla F ( \tilde{y}_{t-1} ) ,  \frac{1}{n} \sum_{j=1}^{n} \nabla f ( y_{j-1}^{(t)} ; \pi^{(t)} ( j ) ) \right \rangle + \frac{\eta_t}{2} \left \| \frac{1}{n} \sum_{j=1}^{n} \nabla f ( y_{j-1}^{(t)} ; \pi^{(t)} ( j ) )  \right \|^2,
\end{align*}
where the last line follows since $\eta_t \leq \frac{1}{L}$. We further have
\begin{align*}
    F(\tilde{x}_t) & \leq (1 - \theta^{(t-1)}) F ( \tilde{x}_{t-1} ) + \theta^{(t-1)} F ( x_{*} )  \\
    & \qquad + \left \langle \nabla F ( \tilde{y}_{t-1} ) - \frac{1}{n} \sum_{j=1}^{n} \nabla f ( y_{j-1}^{(t)} ; \pi^{(t)} ( j ) ) ,  \theta^{(t-1)} ( v^{(t-1)} -  x_{*} ) \right \rangle \\
    & \qquad + \left \langle \frac{1}{n} \sum_{j=1}^{n} \nabla f ( y_{j-1}^{(t)} ; \pi^{(t)} ( j ) ) ,  \theta^{(t-1)} ( v^{(t-1)} -  x_{*} ) \right \rangle \\
    & \qquad -  \eta_t \left \langle \nabla F ( \tilde{y}_{t-1} ) - \frac{1}{n} \sum_{j=1}^{n} \nabla f ( y_{j-1}^{(t)} ; \pi^{(t)} ( j ) )  ,  \frac{1}{n} \sum_{j=1}^{n} \nabla f ( y_{j-1}^{(t)} ; \pi^{(t)} ( j ) ) \right \rangle \\
    & \qquad - \eta_t \left \| \frac{1}{n} \sum_{j=1}^{n} \nabla f ( y_{j-1}^{(t)} ; \pi^{(t)} ( j ) )  \right \|^2 + \frac{\eta_t}{2} \left \| \frac{1}{n} \sum_{j=1}^{n} \nabla f ( y_{j-1}^{(t)} ; \pi^{(t)} ( j ) )  \right \|^2 \\
    & \overset{\eqref{eq_007}}{=} (1 - \theta^{(t-1)}) F ( \tilde{x}_{t-1} ) + \theta^{(t-1)} F ( x_{*} )  \\
    & \qquad + \left \langle \nabla F ( \tilde{y}_{t-1} ) - \frac{1}{n} \sum_{j=1}^{n} \nabla f ( y_{j-1}^{(t)} ; \pi^{(t)} ( j ) ) ,  \theta^{(t-1)} ( v^{(t)} -  x_{*} ) \right \rangle \\
    & \qquad + \left \langle \frac{1}{n} \sum_{j=1}^{n} \nabla f ( y_{j-1}^{(t)} ; \pi^{(t)} ( j ) ) ,  \theta^{(t-1)} ( v^{(t-1)} -  x_{*} ) \right \rangle \\
    & \qquad - \frac{\eta_t}{2} \left \| \frac{1}{n} \sum_{j=1}^{n} \nabla f ( y_{j-1}^{(t)} ; \pi^{(t)} ( j ) )  \right \|^2 \\ 
    & = (1 - \theta^{(t-1)}) F ( \tilde{x}_{t-1} ) + \theta^{(t-1)} F ( x_{*} )  \\
    & \qquad + \left \langle \nabla F ( \tilde{y}_{t-1} ) - \frac{1}{n} \sum_{j=1}^{n} \nabla f ( y_{j-1}^{(t)} ; \pi^{(t)} ( j ) ) ,  \theta^{(t-1)} ( v^{(t)} -  x_{*} ) \right \rangle \\
    & \qquad + \frac{(\theta^{(t-1)})^2}{2 \eta_t} \Big [ \frac{2 \eta_t}{\theta^{(t-1)}} \left \langle \frac{1}{n} \sum_{j=1}^{n} \nabla f ( y_{j-1}^{(t)} ; \pi^{(t)} ( j ) ) ,  ( v^{(t-1)} -  x_{*} ) \right \rangle \\
    & \qquad - \left( \frac{\eta_t}{\theta^{(t-1)}} \right)^2 \left \| \frac{1}{n} \sum_{j=1}^{n} \nabla f ( y_{j-1}^{(t)} ; \pi^{(t)} ( j ) )  \right \|^2 + \| v^{(t-1)} -  x_{*} \|^2 - \| v^{(t-1)} -  x_{*} \|^2 \Big ] \\
    & \overset{\eqref{eq_007}}{=} (1 - \theta^{(t-1)}) F ( \tilde{x}_{t-1} ) + \theta^{(t-1)} F ( x_{*} )  \\
    & \qquad + \left \langle \nabla F ( \tilde{y}_{t-1} ) - \frac{1}{n} \sum_{j=1}^{n} \nabla f ( y_{j-1}^{(t)} ; \pi^{(t)} ( j ) ) ,  \theta^{(t-1)} ( v^{(t)} -  x_{*} ) \right \rangle \\
    & \qquad + \frac{(\theta^{(t-1)})^2}{2 \eta_t} \Big [ \| v^{(t-1)} -  x_{*} \|^2 - \| v^{(t)} -  x_{*} \|^2 \Big ], \tagthis \label{eq_009}
\end{align*}
where we apply equation \eqref{eq_007}.

By the definition $K_t = \frac{1}{n} \sum_{i=1}^{n} \left \|   y_{i}^{(t)} - y_{0}^{(t)}  \right \|^2$, we get that
\begin{align*}
    \left \| \nabla F ( \tilde{y}_{t-1} ) - \frac{1}{n} \sum_{j=1}^{n} \nabla f ( y_{j-1}^{(t)} ; \pi^{(t)} ( j ) ) \right \|^2 & = \left \| \frac{1}{n} \sum_{j=1}^{n} \left ( \nabla f ( \tilde{y}_{t-1} ; \pi^{(t)} ( j ) ) -  \nabla f ( y_{j-1}^{(t)} ; \pi^{(t)} ( j ) ) \right ) \right \|^2 \\
    & \leq \frac{1}{n} \sum_{j=1}^{n} \left \| \nabla f ( \tilde{y}_{t-1} ; \pi^{(t)} ( j ) ) -  \nabla f ( y_{j-1}^{(t)} ; \pi^{(t)} ( j ) )  \right \|^2 \\
    & \overset{\eqref{eq:Lsmooth}}{\leq} L^2 \frac{1}{n} \sum_{j=1}^{n} \left \| y_{0}^{(t)} -  y_{j-1}^{(t)}  \right \|^2 \\
    & \leq \frac{L^2}{n} \sum_{i=1}^{n} \left \|   y_{i}^{(t)} - y_{0}^{(t)}  \right \|^2 = L^2 K_t. \tagthis \label{eq_015}
\end{align*}


From \eqref{eq_009} and using the inequality $\langle a, b \rangle \leq \frac{\| a \|^2}{2 \epsilon_t} + \frac{\epsilon_t \| b \|^2}{2}$ for any $\epsilon_t > 0$, we have
\begin{align*}
     F ( \tilde{x}_{t} ) - F ( x_{*} ) & \leq (1 - \theta^{(t-1)}) [F ( \tilde{x}_{t-1} ) - F ( x_{*} ) ]  \\
    & \qquad + \left \langle \nabla F ( \tilde{y}_{t-1} ) - \frac{1}{n} \sum_{j=1}^{n} \nabla f ( y_{j-1}^{(t)} ; \pi^{(t)} ( j ) ) ,  \theta^{(t-1)} ( v^{(t)} -  x_{*} ) \right \rangle \\
    & \qquad + \frac{(\theta^{(t-1)})^2}{2 \eta_t} \Big [ \| v^{(t-1)} -  x_{*} \|^2 - \| v^{(t)} -  x_{*} \|^2 \Big ] \\
    & \leq (1 - \theta^{(t-1)}) [F ( \tilde{x}_{t-1} ) - F ( x_{*} ) ]  \\
    & \qquad + \frac{\eta_t}{2 \epsilon_t} \left \| \nabla F ( \tilde{y}_{t-1} ) - \frac{1}{n} \sum_{j=1}^{n} \nabla f ( y_{j-1}^{(t)} ; \pi^{(t)} ( j ) ) \right \|^2  + \frac{\epsilon_t (\theta^{(t-1)})^2}{2 \eta_t} \left \|  v^{(t)} -  x_{*}  \right \|^2 \\
    & \qquad + \frac{(\theta^{(t-1)})^2}{2 \eta_t} \Big [ \| v^{(t-1)} -  x_{*} \|^2 - \| v^{(t)} -  x_{*} \|^2 \Big ] \\
    & \overset{\eqref{eq_015}}{\leq} (1 - \theta^{(t-1)}) [F ( \tilde{x}_{t-1} ) - F ( x_{*} ) ]   +  \frac{L^2\eta_t}{2 \epsilon_t}  K_t \\
    & \qquad+ \frac{(\theta^{(t-1)})^2}{2 \eta_t}  \| v^{(t-1)} -  x_{*} \|^2 - \frac{(\theta^{(t-1)})^2}{2 \eta_t} (1-\epsilon_t) \| v^{(t)} -  x_{*} \|^2 .
\end{align*}

Now substituting $\theta^{(t)} = \frac{2}{t+2}$, $\theta^{(t-1)} = \frac{2}{t+1}$ we get

\begin{align*}
     F ( \tilde{x}_{t} ) - F ( x_{*} ) 
    & \leq \frac{t-1}{t+1} [F ( \tilde{x}_{t-1} ) - F ( x_{*} ) ]   +  \frac{L^2\eta_t}{2 \epsilon_t}  K_t \\
    & \qquad+ \frac{2}{\eta_t(t+1)^2}  \| v^{(t-1)} -  x_{*} \|^2 -  \frac{2}{\eta_t(t+1)^2}(1-\epsilon_t) \| v^{(t)} -  x_{*} \|^2. \tagthis \label{eq_016}
\end{align*}

Multiplying two sides by $(t+1)^2$ we have
\begin{align*}
    (t+1)^2 [F ( \tilde{x}_{t} ) - F ( x_{*} ) ]
    & \leq  (t-1)(t+1) [F ( \tilde{x}_{t-1} ) - F ( x_{*} ) ]   +  \frac{L^2\eta_t(t+1)^2}{2 \epsilon_t}  K_t \\
    & \qquad+ \frac{2}{\eta_t}  \| v^{(t-1)} -  x_{*} \|^2 -  \frac{2}{\eta_t}(1-\epsilon_t) \| v^{(t)} -  x_{*} \|^2.
\end{align*}

Summing the previous expression from $t=1$ to $t=T$ we get that
\begin{align*}
    \sum_{t=1}^T(t+1)^2 [F ( \tilde{x}_{t} ) - F ( x_{*} ) ]
    & \leq  \sum_{t=1}^T(t^2-1) [F ( \tilde{x}_{t-1} ) - F ( x_{*} ) ]   +  \sum_{t=1}^T\frac{L^2\eta_t(t+1)^2}{2 \epsilon_t}  K_t \\
    & \qquad+ \sum_{t=1}^T\frac{2}{\eta_t}  \| v^{(t-1)} -  x_{*} \|^2 - \sum_{t=1}^T \frac{2}{\eta_t}(1-\epsilon_t) \| v^{(t)} -  x_{*} \|^2,
\end{align*}
which is equivalent to
\begin{align*}
    \sum_{t=1}^T[(t+1)^2 -1] [F ( \tilde{x}_{t} ) - F ( x_{*} ) ]
    & \leq  \sum_{t=1}^T(t^2-1) [F ( \tilde{x}_{t-1} ) - F ( x_{*} ) ]   +  \sum_{t=1}^T\frac{L^2\eta_t(t+1)^2}{2 \epsilon_t}  K_t \\
    & \qquad+ \sum_{t=1}^T\frac{2}{\eta_t}  \| v^{(t-1)} -  x_{*} \|^2 - \sum_{t=1}^T \frac{2}{\eta_t}(1-\epsilon_t) \| v^{(t)} -  x_{*} \|^2 -\sum_{t=1}^T [F ( \tilde{x}_{t} ) - F ( x_{*} ) ].
\end{align*}
Hence we get the desired estimate of Lemma~\ref{lem_bound_common}:
\begin{align*}
    T(T+2) [F ( \tilde{x}_{T} ) - F ( x_{*} ) ]
    & \leq   \sum_{t=1}^T\frac{L^2\eta_t(t+1)^2}{2 \epsilon_t}  K_t -\sum_{t=1}^T [F ( \tilde{x}_{t} ) - F ( x_{*} ) ]\\
    & \qquad+ \sum_{t=1}^T\frac{2}{\eta_t}  \| v^{(t-1)} -  x_{*} \|^2 - \sum_{t=1}^T \frac{2}{\eta_t}(1-\epsilon_t) \| v^{(t)} -  x_{*} \|^2. \tagthis \label{eq_020}
\end{align*}
\Eproof
\begin{lemma}\label{lem_bound_K_t}
Suppose that Assumption \ref{ass_smooth} holds for \eqref{ERM_problem_01}. 
Let $\sets{\tilde{x}_t}$ be generated by  Algorithm~\ref{shuffling_nesterov_02} with the learning rate $\eta_i^{(t)} := \frac{\eta_t}{n} > 0$ for a given positive sequence $\sets{\eta_t}$ with $\eta_t \leq \frac{1}{2L} $. 
Then we have the following for $t\geq 1$:
\begin{align*}
K_t \leq \frac{8\eta_t^2}{n^3}\sum_{i=1}^{n-1} \left \| \sum_{j=i+1}^n \nabla f ( \tilde{x}_t ; \pi^{(t)} ( j ) ) \right \|^2
    + 4\eta_t^2 \left \|  \nabla F( \tilde{x}_t) \right \|^2.\tagthis \label{eq_017}
\end{align*}
\end{lemma}
\subsection*{Proof of Lemma \ref{lem_bound_K_t}: Bound the term $K_t$}
Let us recall the definition of $K_t$ and $I_t$:
\begin{align*}
K_t &= \frac{1}{n} \sum_{i=1}^{n} \left \|   y_{i}^{(t)} - y_{0}^{(t)}  \right \|^2, \text{ and }
I_t = \frac{1}{n} \sum_{i=1}^{n} \left \|   y_{n}^{(t)} - y_{i}^{(t)}  \right \|^2.
\end{align*}
We consider the individual squared term of $I_t$:
\begin{align*}
    \left \| y_{n}^{(t)} -  y_{i}^{(t)}  \right \|^2 
    &= \frac{\eta_t^2}{n^2}\left \|\sum_{j=i+1}^n \nabla f ( y_{j-1}^{(t)} ; \pi^{(t)} ( j ) )\right \|^2 \\
    &= \frac{\eta_t^2}{n^2}\left \|\sum_{j=i+1}^n \nabla f ( y_{j-1}^{(t)} ; \pi^{(t)} ( j ) ) - \sum_{j=i+1}^n \nabla f ( y_{n}^{(t)} ; \pi^{(t)} ( j ) ) + \sum_{j=i+1}^n \nabla f ( y_{n}^{(t)} ; \pi^{(t)} ( j ) ) \right \|^2 \\
    &\leq \frac{2\eta_t^2}{n^2}\left \|\sum_{j=i+1}^n \nabla f ( y_{j-1}^{(t)} ; \pi^{(t)} ( j ) ) - \sum_{j=i+1}^n \nabla f ( y_{n}^{(t)} ; \pi^{(t)} ( j ) )  \right \|^2 
    + \frac{2\eta_t^2}{n^2}\left \| \sum_{j=i+1}^n \nabla f ( y_{n}^{(t)} ; \pi^{(t)} ( j ) ) \right \|^2 \\
    &\leq \frac{2\eta_t^2}{n^2} (n-i)\sum_{j=i+1}^n \left \| \nabla f ( y_{j-1}^{(t)} ; \pi^{(t)} ( j ) ) -  \nabla f ( y_{n}^{(t)} ; \pi^{(t)} ( j ) )  \right \|^2 
    + \frac{2\eta_t^2}{n^2}\left \| \sum_{j=i+1}^n \nabla f ( y_{n}^{(t)} ; \pi^{(t)} ( j ) ) \right \|^2 ,
\end{align*}
where in the last two lines we use the inequality
 $(u+v)^2 \leq 2u^2 + 2v^2$ and Cauchy-Schwartz inequality. From Assumption \ref{ass_smooth} we have
\begin{align*}
    \left \| y_{n}^{(t)} -  y_{i}^{(t)}  \right \|^2 
     &\leq \frac{2\eta_t^2}{n^2} (n-i)\sum_{j=i+1}^n \left \| \nabla f ( y_{j-1}^{(t)} ; \pi^{(t)} ( j ) ) -  \nabla f ( y_{n}^{(t)} ; \pi^{(t)} ( j ) )  \right \|^2 
    + \frac{2\eta_t^2}{n^2}\left \| \sum_{j=i+1}^n \nabla f ( y_{n}^{(t)} ; \pi^{(t)} ( j ) ) \right \|^2 \\
    &\overset{\eqref{eq:Lsmooth}}{\leq} \frac{2L^2\eta_t^2}{n^2} (n-i)\sum_{j=i+1}^n \left \| y_{j-1}^{(t)} - y_{n}^{(t)}  \right \|^2 
    + \frac{2\eta_t^2}{n^2}\left \| \sum_{j=i+1}^n \nabla f ( y_{n}^{(t)} ; \pi^{(t)} ( j ) ) \right \|^2 \\
    &\leq \frac{2L^2\eta_t^2}{n^2} (n-i) n I_t
    + \frac{2\eta_t^2}{n^2}\left \| \sum_{j=i+1}^n \nabla f ( y_{n}^{(t)} ; \pi^{(t)} ( j ) ) \right \|^2,
\end{align*}
where last inequality follows from definition of $I_t$.
Summing up the previous expression from $i=1$ to $i=n-1$ we get
\begin{align*}
    n I_t = \sum_{i=1}^{n-1} \left \| y_{n}^{(t)} -  y_{i}^{(t)}  \right \|^2 &\leq \frac{2L^2\eta_t^2}{n^2} \sum_{i=1}^{n-1} (n-i) n I_t
    + \frac{2\eta_t^2}{n^2}\sum_{i=1}^{n-1} \left \| \sum_{j=i+1}^n \nabla f ( y_{n}^{(t)} ; \pi^{(t)} ( j ) ) \right \|^2\\
    &\leq \frac{2L^2\eta_t^2}{n^2} \frac{n^2}{2} n I_t
    + \frac{2\eta_t^2}{n^2}\sum_{i=1}^{n-1} \left \| \sum_{j=i+1}^n \nabla f ( y_{n}^{(t)} ; \pi^{(t)} ( j ) ) \right \|^2\\
    &\leq L^2\eta_t^2 n I_t
    + \frac{2\eta_t^2}{n^2}\sum_{i=1}^{n-1} \left \| \sum_{j=i+1}^n \nabla f ( y_{n}^{(t)} ; \pi^{(t)} ( j ) ) \right \|^2,
\end{align*}
where we use the fact that $\sum_{i=1}^{n-1} (n-i)\leq \frac{n^2}{2}$.
Since $\eta_t \leq \frac{1}{2L}$, we have $ \eta_t^2 L^2 \leq \frac{1}{4}$. Hence
\begin{align*}
    \frac{3}{4}n I_t &\leq
    \frac{2\eta_t^2}{n^2}\sum_{i=1}^{n-1} \left \| \sum_{j=i+1}^n \nabla f ( y_{n}^{(t)} ; \pi^{(t)} ( j ) ) \right \|^2,
\end{align*}

and equivalently
\begin{align*}
    I_t &\leq
    \frac{8\eta_t^2}{3n^3}\sum_{i=1}^{n-1} \left \| \sum_{j=i+1}^n \nabla f ( y_{n}^{(t)} ; \pi^{(t)} ( j ) ) \right \|^2. \tagthis \label{eq_021}
\end{align*}

For $i=0$ we have
\begin{align*}
    \left \| y_{n}^{(t)} -  y_{0}^{(t)}  \right \|^2 
    &\leq \frac{2L^2\eta_t^2}{n^2} n^2 I_t
    + \frac{2\eta_t^2}{n^2}\left \| \sum_{j=1}^n \nabla f ( y_{n}^{(t)} ; \pi^{(t)} ( j ) ) \right \|^2\\
    &\overset{\eqref{eq_021}}{\leq} 2L^2\eta_t^2 I_t
    + \frac{2\eta_t^2}{n^2}\left \| \sum_{j=1}^n \nabla f ( y_{n}^{(t)} ; \pi^{(t)} ( j ) ) \right \|^2\\
    &\leq \frac{1}{2} I_t
    + 2\eta_t^2 \left \|  \nabla F( y_{n}^{(t)}) \right \|^2.\tagthis \label{eq_022}
\end{align*}

Now we are ready to investigate $K_t$. By inequality
 $(u+v)^2 \leq 2u^2 + 2v^2$ we get
\begin{align*}
    K_t = \frac{1}{n} \sum_{i=1}^{n} \left \|   y_{i}^{(t)} - y_{0}^{(t)}  \right \|^2
    &\leq \frac{1}{n} \sum_{i=1}^{n} 2\left \|   y_{n}^{(t)} - y_{i}^{(t)}  \right \|^2 + 2\left \|   y_{n}^{(t)} - y_{0}^{(t)}  \right \|^2\\
    &= 2I_t + 2\left \|   y_{n}^{(t)} - y_{0}^{(t)}  \right \|^2\\
    &\overset{\eqref{eq_022}}{\leq} 2I_t + I_t
    + 4\eta_t^2 \left \|  \nabla F( y_{n}^{(t)}) \right \|^2\\
    &\overset{\eqref{eq_021}}{\leq} \frac{8\eta_t^2}{n^3}\sum_{i=1}^{n-1} \left \| \sum_{j=i+1}^n \nabla f ( y_{n}^{(t)} ; \pi^{(t)} ( j ) ) \right \|^2
    + 4\eta_t^2 \left \|  \nabla F( y_{n}^{(t)}) \right \|^2.
\end{align*}

Finally, substituting $y_{n}^{(t)}$ by $\tilde{x}_t$ we get the desired results.
\Eproof

\section{Proof of Theorem \ref{thm_convex_1}: Convex components - Unified schemes}\label{sec_app_unified}
Before proving Theorem \ref{thm_convex_1}, we need the following supplemental Lemma for convex component functions.
\begin{lemma}[Convex component functions]\label{lem_bound_convex_f_i}
Suppose that Assumption \ref{ass_smooth} holds for \eqref{ERM_problem_01} and $f(\cdot;i)$ is convex for every $i \in [n]$. 
Let $\sets{y_i^{(t)}}$ be generated by  Algorithm~\ref{shuffling_nesterov_02} with the learning rate $\eta_i^{(t)} := \frac{\eta_t}{n} > 0$ for a given positive sequence $\sets{\eta_t}$ with $\eta_t \leq \frac{1}{2L} $. 
Then
\begin{align*}
K_t \leq 8\eta_t^2 \left(3L\left(F( \tilde{x}_t  ) - F ( x_*)\right) 
    + \sigma_*^2
     \right). \tagthis \label{eq_019}
\end{align*}
\end{lemma}

\subsection*{Proof of Lemma \ref{lem_bound_convex_f_i}: Bound $K_t$ in terms of the variance $\sigma_*^2$} 
From Lemma \ref{lem_bound_K_t} we have
\begin{align*}
    K_t &\leq \frac{8\eta_t^2}{n^3}\sum_{i=1}^{n-1} \left \| \sum_{j=i+1}^n \nabla f ( \tilde{x}_t ; \pi^{(t)} ( j ) ) \right \|^2
    + 4\eta_t^2 \left \|  \nabla F( \tilde{x}_t) \right \|^2\\
    &= \frac{8\eta_t^2}{n^3}\sum_{i=1}^{n-1} \left \| \sum_{j=i+1}^n \nabla f ( \tilde{x}_t ; \pi^{(t)} ( j ) ) - \sum_{j=i+1}^n \nabla f ( x_* ; \pi^{(t)} ( j ) ) + \sum_{j=i+1}^n \nabla f ( x_* ; \pi^{(t)} ( j ) ) \right \|^2
    + 4\eta_t^2 \left \|  \nabla F( \tilde{x}_t) \right \|^2\\
    &\leq \frac{16\eta_t^2}{n^3}\sum_{i=1}^{n-1} \left \| \sum_{j=i+1}^n \left( \nabla f ( \tilde{x}_t ; \pi^{(t)} ( j ) ) - \nabla f ( x_* ; \pi^{(t)} ( j ) )\right) \right \|^2 
    + \frac{16\eta_t^2}{n^3}\sum_{i=1}^{n-1} \left \| \sum_{j=i+1}^n \nabla f ( x_* ; \pi^{(t)} ( j ) ) \right \|^2
    + 4\eta_t^2 \left \|  \nabla F( \tilde{x}_t) \right \|^2\\
    &\leq \frac{16\eta_t^2}{n^3}\sum_{i=1}^{n-1} (n-i)\sum_{j=i+1}^n \left \|  \nabla f ( \tilde{x}_t ; \pi^{(t)} ( j ) ) - \nabla f ( x_* ; \pi^{(t)} ( j ) )\right \|^2 \\
    & \qquad
    + \frac{16\eta_t^2}{n^3}\sum_{i=1}^{n-1}  (n-i)\sum_{j=i+1}^n\left \| \nabla f ( x_* ; \pi^{(t)} ( j ) ) \right \|^2
    + 4\eta_t^2 \left \|  \nabla F( \tilde{x}_t) \right \|^2,
\end{align*}
where in the last two lines we use the inequality
 $(u+v)^2 \leq 2u^2 + 2v^2$ and Cauchy-Schwartz inequality. By the definition of $D_t$ we have
\begin{align*}
    K_t &\leq \frac{16\eta_t^2}{n^3}\sum_{i=1}^{n-1} (n-i) D_t 
    + \frac{16\eta_t^2}{n^3}\sum_{i=1}^{n-1}  (n-i)\sum_{j=1}^n\left \| \nabla f ( x_* ; \pi^{(t)} ( j ) ) \right \|^2
    + 4\eta_t^2 \left \|  \nabla F( \tilde{x}_t) \right \|^2\\
    &\leq \frac{8\eta_t^2}{n} D_t 
    + \frac{8\eta_t^2}{n} n \sigma_*^2
    + 4\eta_t^2 \left \|  \nabla F( \tilde{x}_t) \right \|^2,
\end{align*}
where we use the fact that $\sum_{i=1}^{n-1} (n-i)\leq \frac{n^2}{2}$.

Let us consider the term $D_t$. Since $f_i$ is convex, we have the following for every $t \geq 1$
\begin{align*}
    D_t &= \sum_{j=1}^{n}\left\|  \nabla f ( \tilde{x}_t ; \pi^{(t)} ( j ) ) - \nabla f ( x_*; \pi^{(t)} ( j ) ) ) \right\|^2\\
    &\leq 2L \sum_{j=1}^{n} \left( f ( \tilde{x}_t ; \pi^{(t)} ( j ) ) - f ( x_*; \pi^{(t)} ( j ) ) - \iprods{\nabla f ( x_*; \pi^{(t)} ( j ) ), \tilde{x}_t - x_*}\right)\\
    &\leq 2nL \left( F( \tilde{x}_t  ) - F ( x_*) - \iprods{\nabla F ( x_* ), \tilde{x}_t - x_*}\right)\\
    &= 2nL \left(F( \tilde{x}_t  ) - F ( x_*)\right).
\end{align*}
Substitute this to the previous equation we get:
\begin{align*}
    K_t &\leq \frac{8\eta_t^2}{n} D_t 
    + \frac{8\eta_t^2}{n} n \sigma_*^2
    + 4\eta_t^2 \left \|  \nabla F( \tilde{x}_t) \right \|^2\\
    &\leq 16L\eta_t^2\left(F( \tilde{x}_t  ) - F ( x_*)\right) 
    + 8\eta_t^2\sigma_*^2
    + 4\eta_t^2 \left \|  \nabla F( \tilde{x}_t) \right \|^2.
\end{align*}
Since $F$ is $L$-smoooth and convex, we have $\left \|  \nabla F( \tilde{x}_t) \right \|^2 \leq 2L\left(F( \tilde{x}_t  ) - F (x_*)\right)$ \cite{Nesterov2004}. Hence
\begin{align*}
    K_t 
    &\leq 16L\eta_t^2\left(F( \tilde{x}_t  ) - F ( x_*)\right) 
    + 8\eta_t^2\sigma_*^2
    + 4\eta_t^2 \cdot 2L\left(F( \tilde{x}_t  ) - F (x_*)\right)\\
    &\leq 24L\eta_t^2\left(F( \tilde{x}_t  ) - F ( x_*)\right) 
    + 8\eta_t^2\sigma_*^2.
\end{align*}
Thus we have the estimate of Lemma \ref{lem_bound_convex_f_i}.
\Eproof

\subsection*{Proof of Theorem \ref{thm_convex_1}} 

Let us start with inequality \eqref{eq_020_0} from Lemma \ref{lem_bound_common}. Applying Lemma \ref{lem_bound_convex_f_i} we have
\begin{align*}
    T(T+2) [F ( \tilde{x}_{T} ) - F ( x_{*} ) ]
    & \leq   \sum_{t=1}^T\frac{L^2\eta_t(t+1)^2}{2 \epsilon_t}  K_t -\sum_{t=1}^T [F ( \tilde{x}_{t} ) - F ( x_{*} ) ]\\
    & \qquad+ \sum_{t=1}^T\frac{2}{\eta_t}  \| v^{(t-1)} -  x_{*} \|^2 - \sum_{t=1}^T \frac{2}{\eta_t}(1-\epsilon_t) \| v^{(t)} -  x_{*} \|^2 \\
    & \overset{\eqref{eq_019}}{\leq}    \sum_{t=1}^T\frac{4L^2\eta_t^3(t+1)^2}{\epsilon_t}  \left(3L\left(F( \tilde{x}_t  ) - F ( x_*)\right) 
    + \sigma_*^2 \right) -\sum_{t=1}^T [F ( \tilde{x}_{t} ) - F ( x_{*} ) ]\\
    & \qquad+ \sum_{t=1}^T\frac{2}{\eta_t}  \| v^{(t-1)} -  x_{*} \|^2 - \sum_{t=1}^T \frac{2}{\eta_t}(1-\epsilon_t) \| v^{(t)} -  x_{*} \|^2.
\end{align*}
From the choice $\eta_t = \frac{k\alpha^t}{LT}$ we have $\frac{2}{\eta_t} = \frac{2LT}{k\alpha^t}$ and
\begin{align*}
     T(T+2) [F ( \tilde{x}_{T} ) - F ( x_{*} ) ]&\leq   \sum_{t=1}^T\frac{k^3\alpha^{3t}}{L^3T^3} \frac{4L^2(t+1)^2}{\epsilon_t}  \left(3L\left(F( \tilde{x}_t  ) - F ( x_*)\right) + \sigma_*^2  \right) -\sum_{t=1}^T [F ( \tilde{x}_{t} ) - F ( x_{*} ) ]\\
    & \quad+ \sum_{t=1}^T \frac{2LT}{k\alpha^t}  \| v^{(t-1)} -  x_{*} \|^2 - \sum_{t=1}^T  \frac{2LT}{k\alpha^t}(1-\epsilon_t) \| v^{(t)} -  x_{*} \|^2.
\end{align*}
In addition, we choose $\epsilon_t = \frac{\alpha -1}{\alpha}$ and $(1-\epsilon_t) = \frac{1}{\alpha}$. The last two terms cancel out that
\begin{align*}
     T(T+2) [F ( \tilde{x}_{T} ) - F ( x_{*} ) ]
    &\leq   \sum_{t=1}^T\frac{k^3\alpha^{3t}}{L^3T^3} \frac{4\alpha L^2(t+1)^2}{\alpha -1}  \left(3L\left(F( \tilde{x}_t  ) - F ( x_*)\right) + \sigma_*^2  \right) -\sum_{t=1}^T [F ( \tilde{x}_{t} ) - F ( x_{*} ) ]\\
    & \quad+ \sum_{t=1}^T \frac{2LT}{k\alpha^t}  \| v^{(t-1)} -  x_{*} \|^2 - \sum_{t=1}^T  \frac{2LT}{k\alpha^{t+1}}\| v^{(t)} -  x_{*} \|^2.\\
    &\leq   \sum_{t=1}^T\frac{k^3\alpha^{3t+1}}{LT^3} \frac{4(t+1)^2}{\alpha-1}  \left(3L\left(F( \tilde{x}_t  ) - F (x_*)\right) + \sigma_*^2 \right) -\sum_{t=1}^T [F ( \tilde{x}_{t} ) - F ( x_{*} ) ]\\
    & \qquad+ \frac{2LT}{k\alpha}  \| v^{0} -  x_{*} \|^2.
\end{align*}
Note that $\alpha = 1+ \frac{1}{T}$ ($1 \leq \alpha \leq \frac{3}{2}$ for $T \geq 2$). Hence $\alpha -1= \frac{1}{T}$, $\alpha^t \leq \alpha^T = \left( 1+ \frac{1}{T}\right)^T \leq e$ and
\begin{align*}
     T(T+2) [F ( \tilde{x}_{T} ) - F ( x_{*} ) ]
    &\leq   \sum_{t=1}^T\frac{k^3e^3\alpha}{LT^2}  4 (t+1)^2 \left(3L\left(F( \tilde{x}_t  ) - F (x_*)\right) + \sigma_*^2  \right) -\sum_{t=1}^T [F ( \tilde{x}_{t} ) - F ( x_{*} ) ]
    + \frac{2LT}{k\alpha}  \| v^{0} -  x_{*} \|^2 \\
    &\leq   \sum_{t=1}^T \left[\frac{12k^3e^3\alpha(t+1)^2}{T^2}    -1\right] [F ( \tilde{x}_{t} ) - F ( x_{*} ) ] + \sum_{t=1}^T\frac{4k^3e^3\alpha (t+1)^2  \sigma_*^2 }{LT^2}   
    + \frac{2LT}{k\alpha}  \| v^{0} -  x_{*} \|^2 .
\end{align*}

From the choice $k = \frac{1}{e \alpha \sqrt[3]{12}}$, we have $12k^3 e^3 \alpha^3 = 1$.
Hence for every $t \geq 1$ we have 
\begin{align*}
    \frac{12k^3e^3\alpha(t+1)^2}{T^2} -1 \leq \frac{12k^3e^3\alpha(T+1)^2}{T^2} -1  \leq 12k^3e^3 \alpha^3 -1 = 0,
\end{align*}
where we use the fact that $\alpha = 1+ \frac{1}{T} =  \frac{T+1}{T}$.

We further have 
\begin{align*}
     T(T+2) [F ( \tilde{x}_{T} ) - F ( x_{*} ) ]
    &\leq  \sum_{t=1}^T\frac{4k^3e^3 \alpha  (t+1)^2  \sigma_*^2 }{LT^2}   
    + \frac{2LT}{k\alpha}  \| v^{0} -  x_{*} \|^2 \\
    &\leq \frac{4k^3e^3 \alpha  (T+2)^3  \sigma_*^2 }{3LT^2}   
    + \frac{2LT}{k\alpha}  \| v^{0} -  x_{*} \|^2,
\end{align*}
where we use the fact that $\sum_{t=1}^{T} (t+1)^2 \leq \frac{(T+2)^3}{3}$. Dividing both sides by $T(T+2)$ and substituting $k = \frac{1}{e \alpha \sqrt[3]{12}}$ and $12k^3 e^3 \alpha^3 = 1$ we have
\begin{align*}
    F ( \tilde{x}_{T} ) - F ( x_{*} ) 
    &\leq \frac{4k^3e^3 \alpha  (T+2)^2  \sigma_*^2 }{3LT^3}   
    + \frac{2L}{k\alpha(T+2)}  \| v^{0} -  x_{*} \|^2\\
    &\leq \frac{ (T+2)^2  \sigma_*^2 }{9 \alpha^2 LT^3}   
    + \frac{2Le \sqrt[3]{12}}{T+2}  \| v^{0} -  x_{*} \|^2\\
    &\leq \frac{4 \sigma_*^2 }{9  LT}   
    + \frac{2Le \sqrt[3]{12}}{T}  \| v^{0} -  x_{*} \|^2, 
\end{align*}
where $(T+2)^2 \leq 4T^2$ for $T \geq 2$. Note that $v^{0} = \tilde{x}_0$, we get the desired results.
\Eproof
\subsection*{Proof of Corollary \ref{cor_comp_unified}: Computational complexity of Theorem \ref{thm_convex_1}} 
\begin{corollary}\label{cor_comp_unified}
Assume the same conditions as in Theorem \ref{thm_convex_1}, i.e. Assumption \ref{ass_smooth} and \ref{ass_convex} holds for \eqref{ERM_problem_01}. The computational complexity needed by Algorithm~\ref{shuffling_nesterov_02} to reach an $\epsilon$-accurate solution $x$ that satisfies $F(x) - F(x_*) \leq \epsilon$ is
\begin{align}
    nT
    = \Ocal \left( \frac{n\sigma_*^2}{L\epsilon}+ \frac{nL\| \tilde{x}_0 -  x_{*} \|^2}{\epsilon}\right).
    \label{eq_comp_unified}
\end{align} 
\end{corollary}
By Theorem \ref{thm_convex_1} we have 
\begin{align*}
    F ( \tilde{x}_{T} ) - F ( x_{*} ) 
    &\leq \frac{4 \sigma_*^2 }{9  LT}   
    + \frac{2Le \sqrt[3]{12}}{T}  \| \tilde{x}_0 -  x_{*} \|^2.
\end{align*} 
In order to reach an $\epsilon$-accurate solution $x = \tilde{x}_{T}$ that satisfies $F(x) - F(x_*) \leq \epsilon$, we need
\begin{align*}
    \frac{4 \sigma_*^2 }{9  LT}   
     &\leq \frac{\epsilon}{2} \text{ and  }
    \frac{2Le \sqrt[3]{12}}{T}  \| \tilde{x}_0 -  x_{*} \|^2 \leq \frac{\epsilon}{2},
\end{align*} 
which is equivalent to 
\begin{align*}
    T
     &\geq \frac{8\sigma_*^2}{9L\epsilon} \text{ and  }
    T  \geq \frac{4Le \sqrt[3]{12}\| \tilde{x}_0 -  x_{*} \|^2}{\epsilon}.
\end{align*} 
Hence the number of individual gradient evaluations needed is
\begin{align*}
    nT
     = \max \left( \frac{8n\sigma_*^2}{9L\epsilon}, \frac{4nLe \sqrt[3]{12}\| \tilde{x}_0 -  x_{*} \|^2}{\epsilon}\right)\leq  \frac{8n\sigma_*^2}{9L\epsilon}+ \frac{4nLe \sqrt[3]{12}\| \tilde{x}_0 -  x_{*} \|^2}{\epsilon}= \Ocal \left( \frac{n\sigma_*^2}{L\epsilon}+ \frac{nL\| \tilde{x}_0 -  x_{*} \|^2}{\epsilon}\right).
\end{align*} 
\Eproof







\section{Proof of Theorem \ref{thm_convex_2}: Bounded variance - Unified schemes}\label{sec_app_variance}
Before proving Theorem \ref{thm_convex_2}, we need the following supplemental Lemma for generalized bounded variance assumption. 
\begin{lemma}[Bounded variance]\label{lem_bound_variance}
Suppose that Assumption \ref{ass_smooth} and \ref{ass_bounded_variance} holds for \eqref{ERM_problem_01}. 
Let $\sets{\tilde{x}_t}$ be generated by  Algorithm~\ref{shuffling_nesterov_02} with the learning rate $\eta_i^{(t)} := \frac{\eta_t}{n} > 0$ for a given positive sequence $\sets{\eta_t}$ with $\eta_t \leq \frac{1}{2L} $. 
Then
\begin{align*}
K_t \leq \frac{4\eta_t^2}{3}\left((6\Theta + 7) \| \nabla F(\tilde{x}_t) \|^2 + 6\sigma^2\right).\tagthis \label{eq_018}
\end{align*}
\end{lemma}
\subsection*{Proof of Lemma \ref{lem_bound_variance}: Bound $K_t$ in terms of the variance $\sigma^2$} 
From Lemma \ref{lem_bound_K_t} we have
\begin{align*}
    K_t &\leq \frac{8\eta_t^2}{n^3}\sum_{i=1}^{n-1} \left \| \sum_{j=i+1}^n \nabla f ( \tilde{x}_t ; \pi^{(t)} ( j ) ) \right \|^2
    + 4\eta_t^2 \left \|  \nabla F( \tilde{x}_t) \right \|^2 \\
    &= \frac{8\eta_t^2}{n^3}\sum_{i=1}^{n-1} \left \| \sum_{j=i+1}^n \left(\nabla f ( \tilde{x}_t ; \pi^{(t)} ( j ) ) - \nabla F( \tilde{x}_t) + \nabla F( \tilde{x}_t)\right )\right \|^2
    + 4\eta_t^2 \left \|  \nabla F( \tilde{x}_t) \right \|^2 \\
    &\leq \frac{16\eta_t^2}{n^3}\sum_{i=1}^{n-1} \left \| \sum_{j=i+1}^n \left(\nabla f ( \tilde{x}_t ; \pi^{(t)} ( j ) ) - \nabla F( \tilde{x}_t)\right )\right \|^2 + \frac{16\eta_t^2}{n^3}\sum_{i=1}^{n-1} \left \| \sum_{j=i+1}^n\nabla F( \tilde{x}_t) \right \|^2
    + 4\eta_t^2 \left \|  \nabla F( \tilde{x}_t) \right \|^2 \\
    &\leq \frac{16\eta_t^2}{n^3}\sum_{i=1}^{n-1} (n-i) \sum_{j=i+1}^n\left \|  \nabla f ( \tilde{x}_t ; \pi^{(t)} ( j ) ) - \nabla F( \tilde{x}_t)\right \|^2 + \frac{16\eta_t^2}{n^3}\sum_{i=1}^{n-1}(n-i)^2 \left \| \nabla F( \tilde{x}_t) \right \|^2
    + 4\eta_t^2 \left \|  \nabla F( \tilde{x}_t) \right \|^2 ,
\end{align*}
where in the last two lines we use the inequality
 $(u+v)^2 \leq 2u^2 + 2v^2$ and Cauchy-Schwartz inequality. By Assumption \ref{ass_bounded_variance} we have
\begin{align*}
    K_t &\leq \frac{16\eta_t^2}{n^3}\sum_{i=1}^{n-1} (n-i) n\left(\Theta \| \nabla F(\tilde{x}_t) \|^2 + \sigma^2\right) + \frac{16\eta_t^2}{n^3}\sum_{i=1}^{n-1}(n-i)^2 \left \| \nabla F( \tilde{x}_t) \right \|^2
    + 4\eta_t^2 \left \|  \nabla F( \tilde{x}_t) \right \|^2 \\
    &\leq 8\eta_t^2 \left(\Theta \| \nabla F(\tilde{x}_t) \|^2 + \sigma^2\right) + \frac{16\eta_t^2}{3} \left \| \nabla F( \tilde{x}_t) \right \|^2
    + 4\eta_t^2 \left \|  \nabla F( \tilde{x}_t) \right \|^2 \\
    &\leq \frac{4\eta_t^2}{3}\left((6\Theta+ 7) \| \nabla F(\tilde{x}_t) \|^2 + 6\sigma^2\right),
\end{align*}
where we use the inequalities  $\sum_{i=1}^{n-1} (n-i)\leq \frac{n^2}{2}$
and $\sum_{i=1}^{n-1} (n-i)^2\leq \frac{n^3}{3}$.

\subsection*{Proof of Theorem \ref{thm_convex_2}} 
Let us start with inequality \eqref{eq_020_0} from Lemma \ref{lem_bound_common}. Applying Lemma \ref{lem_bound_variance} we have
\begin{align*}
    T(T+2) [F ( \tilde{x}_{T} ) - F ( x_{*} ) ]
    & \leq   \sum_{t=1}^T\frac{L^2\eta_t(t+1)^2}{2 \epsilon_t}  K_t -\sum_{t=1}^T [F ( \tilde{x}_{t} ) - F ( x_{*} ) ]\\
    & \qquad+ \sum_{t=1}^T\frac{2}{\eta_t}  \| v^{(t-1)} -  x_{*} \|^2 - \sum_{t=1}^T \frac{2}{\eta_t}(1-\epsilon_t) \| v^{(t)} -  x_{*} \|^2 \\
    & \leq   \sum_{t=1}^T\frac{L^2\eta_t(t+1)^2}{2 \epsilon_t} \frac{4\eta_t^2}{3}\left((6\Theta + 7) \| \nabla F(\tilde{x}_t) \|^2 + 6\sigma^2\right) -\sum_{t=1}^T [F ( \tilde{x}_{t} ) - F ( x_{*} ) ]\\
    & \qquad+ \sum_{t=1}^T\frac{2}{\eta_t}  \| v^{(t-1)} -  x_{*} \|^2 - \sum_{t=1}^T \frac{2}{\eta_t}(1-\epsilon_t) \| v^{(t)} -  x_{*} \|^2 \\
    & \leq   \sum_{t=1}^T\frac{2L^2\eta_t^3(t+1)^2}{3\epsilon_t} \left((6\Theta + 7) \| \nabla F(\tilde{x}_t) \|^2 + 6\sigma^2\right) -\frac{1}{2L}\sum_{t=1}^T \| \nabla F(\tilde{x}_t) \|^2\\
    & \qquad+ \sum_{t=1}^T\frac{2}{\eta_t}  \| v^{(t-1)} -  x_{*} \|^2 - \sum_{t=1}^T \frac{2}{\eta_t}(1-\epsilon_t) \| v^{(t)} -  x_{*} \|^2,
\end{align*}
where we use the inequality $F ( \tilde{x}_{t} ) - F ( x_{*} ) \geq \frac{1}{2L}\| \nabla F(\tilde{x}_t) \|^2$ since $F$ is $L$-smoooth and convex \cite{Nesterov2004}.

From the choice $\eta_t = \frac{k\alpha^t}{LT}$ we have $\frac{2}{\eta_t} = \frac{2LT}{k\alpha^t}$ and
\begin{align*}
    T(T+2) [F ( \tilde{x}_{T} ) - F ( x_{*} ) ]
    & \leq   \sum_{t=1}^T\frac{k^3\alpha^{3t}}{L^3T^3} \frac{2L^2(t+1)^2}{3\epsilon_t} \left((6\Theta + 7) \| \nabla F(\tilde{x}_t) \|^2 + 6\sigma^2\right) -\frac{1}{2L}\sum_{t=1}^T \| \nabla F(\tilde{x}_t) \|^2\\
    & \qquad+ \sum_{t=1}^T\frac{2LT}{k\alpha^t}  \| v^{(t-1)} -  x_{*} \|^2 - \sum_{t=1}^T \frac{2LT}{k\alpha^t}(1-\epsilon_t) \| v^{(t)} -  x_{*} \|^2.
\end{align*}
In addition, we choose $\epsilon_t = \frac{\alpha -1}{\alpha}$ and $(1-\epsilon_t) = \frac{1}{\alpha}$. The last two terms cancel out that
\begin{align*}
    T(T+2) [F ( \tilde{x}_{T} ) - F ( x_{*} ) ]
    & \leq   \sum_{t=1}^T\frac{k^3\alpha^{3t}}{L^3T^3} \frac{2\alpha L^2(t+1)^2}{3(\alpha-1)} \left((6\Theta + 7) \| \nabla F(\tilde{x}_t) \|^2 + 6\sigma^2\right) -\frac{1}{2L}\sum_{t=1}^T \| \nabla F(\tilde{x}_t) \|^2\\
    & \qquad+ \sum_{t=1}^T\frac{2LT}{k\alpha^t}  \| v^{(t-1)} -  x_{*} \|^2 - \sum_{t=1}^T \frac{2LT}{k\alpha^{t+1}} \| v^{(t)} -  x_{*} \|^2\\
    & \leq   \sum_{t=1}^T\frac{k^3\alpha^{3t+1}}{LT^3} \frac{2(t+1)^2}{3(\alpha-1)} \left((6\Theta + 7) \| \nabla F(\tilde{x}_t) \|^2 + 6\sigma^2\right) -\frac{1}{2L}\sum_{t=1}^T \| \nabla F(\tilde{x}_t) \|^2 \\ & \qquad + \frac{2LT}{k\alpha}  \| v^{0} -  x_{*} \|^2 .
\end{align*}

Note that $\alpha = 1+ \frac{1}{T}$ ($1 \leq \alpha \leq \frac{3}{2}$ for $T \geq 2$). Hence $\alpha -1= \frac{1}{T}$, $\alpha^t \leq \alpha^T = \left( 1+ \frac{1}{T}\right)^T \leq e$ and
\begin{align*}
    T(T+2) [F ( \tilde{x}_{T} ) - F ( x_{*} ) ]
    & \leq   \sum_{t=1}^T\frac{2k^3e^3\alpha(t+1)^2}{3LT^2}  \left((6\Theta + 7) \| \nabla F(\tilde{x}_t) \|^2 + 6\sigma^2\right) -\frac{1}{2L}\sum_{t=1}^T \| \nabla F(\tilde{x}_t) \|^2 \\ & \qquad + \frac{2LT}{k\alpha}  \| v^{0} -  x_{*} \|^2 \\
    & \leq   \sum_{t=1}^T \left[\frac{2k^3e^3\alpha(t+1)^2(6\Theta + 7) }{3LT^2}   -\frac{1}{2L} \right]\| \nabla F(\tilde{x}_t) \|^2\\
    & \qquad + \sum_{t=1}^T\frac{4k^3e^3\alpha(t+1)^2}{LT^2} \sigma^2  + \frac{2LT}{k\alpha}  \| v^{0} -  x_{*} \|^2 .
\end{align*}
From the choice $k = \frac{1}{e \alpha \sqrt[3]{2(6\Theta + 7)}}$, we have $2k^3e^3\alpha^3(6\Theta + 7) = 1 $. Hence for every $t \geq 1$ we have 
\begin{align*}
    \frac{2k^3e^3\alpha(t+1)^2(6\Theta + 7) }{3LT^2} - \frac{1}{2L} \leq
    \frac{2k^3e^3\alpha(T+1)^2(6\Theta + 7) }{3LT^2} - \frac{1}{2L}\leq
    \frac{2k^3e^3\alpha^3(6\Theta + 7) }{3L} - \frac{1}{2L} \leq
    \frac{1 }{3L}  - \frac{1}{2L} \leq 0.
\end{align*}
where we use the fact that $\alpha = 1+ \frac{1}{T} =  \frac{T+1}{T}$.

We further have
\begin{align*}
    T(T+2) [F ( \tilde{x}_{T} ) - F ( x_{*} ) ]
    & \leq   \sum_{t=1}^T\frac{4k^3e^3\alpha(t+1)^2}{LT^2} \sigma^2  + \frac{2LT}{k\alpha}  \| v^{0} -  x_{*} \|^2 \\
    & \leq   \frac{4k^3e^3\alpha(T+2)^3}{3LT^2} \sigma^2  + \frac{2LT}{k\alpha}  \| v^{0} -  x_{*} \|^2 ,
\end{align*}
where we use the fact that $\sum_{t=1}^{T} (t+1)^2 \leq \frac{(T+2)^3}{3}$. Dividing both sides by $T(T+2)$ and substituting $k = \frac{1}{e \alpha \sqrt[3]{2(6\Theta + 7)}}$ and $2k^3e^3\alpha^3(6\Theta + 7) = 1 $ we have
\begin{align*}
    F ( \tilde{x}_{T} ) - F ( x_{*} ) 
    & \leq   \frac{4k^3e^3\alpha(T+2)^2}{3LT^3} \sigma^2  + \frac{2L}{k\alpha(T+2)}  \| v^{0} -  x_{*} \|^2\\
    & \leq   \frac{2(T+2)^2\sigma^2 }{3\alpha^2 (6\Theta + 7)LT^3}   + \frac{2Le  \sqrt[3]{2(6\Theta + 7)} }{T+2}  \| v^{0} -  x_{*} \|^2\\
    & \leq   \frac{8\sigma^2 }{3(6\Theta + 7)LT}   + \frac{2Le  \sqrt[3]{2(6\Theta + 7)} }{T}  \| v^{0} -  x_{*} \|^2,
\end{align*}
where $(T+2)^2 \leq 4T^2$ for $T \geq 2$. Note that $v^{0} = \tilde{x}_0$, we get the desired results.
\Eproof

\section{Proof of Theorem \ref{thm_convex_RR}: Convex components - Randomized schemes}\label{sec_app_RR} 
Before proving Theorem \ref{thm_convex_RR}, we need two supplemental Lemmas for Randomized sampling schemes.
The first Lemma is \citep{mishchenko2020random}[Lemma 1]  for sampling without replacement.

\begin{lemma}[Lemma 1 in \cite{mishchenko2020random}]
Let $X_1, \cdots, X_n \in \R^d$ be fixed vectors, $\bar{X} := \frac{1}{n} \sum_{i=1}^n X_i$ be their average and $\sigma^2 := \frac{1}{n} \sum_{i=1}^n \|X_i -\bar{X}\|^2$ 
be the population variance. Fix any $k \in \{1,\cdots, n\}$, let $X_{\pi_1}, \cdots, X_{\pi_k}$ be sampled uniformly without replacement from $\{X_1, \cdots, X_n\}$ and $\bar{X}_\pi$ be their average. 
Then, the sample average and the variance are given, respectively by
\begin{align*}
    \E [\bar{X}_\pi] = \bar{X} \qquad \text{and} \qquad \E \left[ \|\bar{X}_\pi - \bar{X}\|^2 \right] = \frac{n-k}{k(n-1)} \sigma^2.
\end{align*}
\end{lemma}
Using this result we are able to prove the next Lemma~\ref{lem_bound_convex_RR} as follows.
\begin{lemma}[Randomized Sampling]\label{lem_bound_convex_RR}
Suppose that Assumption \ref{ass_smooth} holds for \eqref{ERM_problem_01} and $f(\cdot;i)$ is convex for every $i \in [n]$. 
Let $\sets{y_i^{(t)}}$ be generated by  Algorithm~\ref{shuffling_nesterov_02} with the learning rate $\eta_i^{(t)} := \frac{\eta_t}{n} > 0$ for a given positive sequence $\sets{\eta_t}$ with $\eta_t \leq \frac{1}{2L} $. 
Then
\begin{align*}
\E [K_t] \leq 8\eta_t^2 \left(3L \E \left[F( \tilde{x}_t  ) - F ( x_*)\right]
    + \frac{2\sigma_{*}^2}{3n}
     \right). \tagthis \label{eq_019_RR}
\end{align*}
\end{lemma}

\subsection*{Proof of Lemma \ref{lem_bound_convex_f_i}: Bound $K_t$ in terms of the variance $\sigma_*^2$} 
From Lemma \ref{lem_bound_K_t} we have
\begin{align*}
    K_t &\leq \frac{8\eta_t^2}{n^3}\sum_{i=1}^{n-1} \left \| \sum_{j=i+1}^n \nabla f ( \tilde{x}_t ; \pi^{(t)} ( j ) ) \right \|^2
    + 4\eta_t^2 \left \|  \nabla F( \tilde{x}_t) \right \|^2\\
    &= \frac{8\eta_t^2}{n^3}\sum_{i=1}^{n-1} \left \| \sum_{j=i+1}^n \nabla f ( \tilde{x}_t ; \pi^{(t)} ( j ) ) - \sum_{j=i+1}^n \nabla f ( x_* ; \pi^{(t)} ( j ) ) + \sum_{j=i+1}^n \nabla f ( x_* ; \pi^{(t)} ( j ) ) \right \|^2
    + 4\eta_t^2 \left \|  \nabla F( \tilde{x}_t) \right \|^2\\
    &\leq \frac{16\eta_t^2}{n^3}\sum_{i=1}^{n-1} \left \| \sum_{j=i+1}^n \left( \nabla f ( \tilde{x}_t ; \pi^{(t)} ( j ) ) - \nabla f ( x_* ; \pi^{(t)} ( j ) )\right) \right \|^2 
    + \frac{16\eta_t^2}{n^3}\sum_{i=1}^{n-1} \left \| \sum_{j=i+1}^n \nabla f ( x_* ; \pi^{(t)} ( j ) ) \right \|^2
    + 4\eta_t^2 \left \|  \nabla F( \tilde{x}_t) \right \|^2\\
    &\leq \frac{16\eta_t^2}{n^3}\sum_{i=1}^{n-1} (n-i)\sum_{j=i+1}^n \left \|  \nabla f ( \tilde{x}_t ; \pi^{(t)} ( j ) ) - \nabla f ( x_* ; \pi^{(t)} ( j ) )\right \|^2 \\
    & \qquad
    + \frac{16\eta_t^2}{n^3}\sum_{i=1}^{n-1} \left \| \sum_{j=i+1}^n \nabla f ( x_* ; \pi^{(t)} ( j ) ) \right \|^2
    + 4\eta_t^2 \left \|  \nabla F( \tilde{x}_t) \right \|^2,
\end{align*}
where in the last two lines we use the inequality
 $(u+v)^2 \leq 2u^2 + 2v^2$ and Cauchy-Schwartz inequality. By the definition of $D_t$ we have
\begin{align*}
    K_t &\leq \frac{16\eta_t^2}{n^3}\sum_{i=1}^{n-1} (n-i) D_t 
    + \frac{16\eta_t^2}{n^3}\sum_{i=1}^{n-1} \left \| \sum_{j=i+1}^n \nabla f ( x_* ; \pi^{(t)} ( j ) ) \right \|^2
    + 4\eta_t^2 \left \|  \nabla F( \tilde{x}_t) \right \|^2\\
    &\leq \frac{8\eta_t^2}{n} D_t 
    + \frac{16\eta_t^2}{n^3}\sum_{i=1}^{n-1} \left \| \sum_{j=i+1}^n \nabla f ( x_* ; \pi^{(t)} ( j ) ) \right \|^2
    + 4\eta_t^2 \left \|  \nabla F( \tilde{x}_t) \right \|^2,
\end{align*}
where we use the fact that $\sum_{i=1}^{n-1} (n-i)\leq \frac{n^2}{2}$.

Let us consider the term $D_t$. Since $f_i$ is convex, we have the following for every $t \geq 1$
\begin{align*}
    D_t &= \sum_{j=1}^{n}\left\|  \nabla f ( \tilde{x}_t ; \pi^{(t)} ( j ) ) - \nabla f ( x_*; \pi^{(t)} ( j ) ) ) \right\|^2\\
    &\leq 2L \sum_{j=1}^{n} \left( f ( \tilde{x}_t ; \pi^{(t)} ( j ) ) - f ( x_*; \pi^{(t)} ( j ) ) - \iprods{\nabla f ( x_*; \pi^{(t)} ( j ) ), \tilde{x}_t - x_*}\right)\\
    &\leq 2nL \left( F( \tilde{x}_t  ) - F ( x_*) - \iprods{\nabla F ( x_* ), \tilde{x}_t - x_*}\right)\\
    &= 2nL \left(F( \tilde{x}_t  ) - F ( x_*)\right).
\end{align*}
Substitute this to the previous equation we get:
\begin{align*}
    K_t &\leq \frac{8\eta_t^2}{n} D_t 
    + \frac{16\eta_t^2}{n^3}\sum_{i=1}^{n-1} \left \| \sum_{j=i+1}^n \nabla f ( x_* ; \pi^{(t)} ( j ) ) \right \|^2
    + 4\eta_t^2 \left \|  \nabla F( \tilde{x}_t) \right \|^2\\
    &\leq 16L\eta_t^2\left(F( \tilde{x}_t  ) - F ( x_*)\right) 
    + \frac{16\eta_t^2}{n^3}\sum_{i=1}^{n-1} \left \| \sum_{j=i+1}^n \nabla f ( x_* ; \pi^{(t)} ( j ) ) \right \|^2
    + 4\eta_t^2 \left \|  \nabla F( \tilde{x}_t) \right \|^2.
\end{align*}
Since $F$ is $L$-smoooth and convex, we have $\left \|  \nabla F( \tilde{x}_t) \right \|^2 \leq 2L\left(F( \tilde{x}_t  ) - F (x_*)\right)$ \cite{Nesterov2004}. Hence
\begin{align*}
    K_t 
    &\leq 16L\eta_t^2\left(F( \tilde{x}_t  ) - F ( x_*)\right) 
    + \frac{16\eta_t^2}{n^3}\sum_{i=1}^{n-1} \left \| \sum_{j=i+1}^n \nabla f ( x_* ; \pi^{(t)} ( j ) ) \right \|^2
    + 4\eta_t^2 \cdot 2L\left(F( \tilde{x}_t  ) - F (x_*)\right)\\
    &\leq 24L\eta_t^2\left(F( \tilde{x}_t  ) - F ( x_*)\right) 
    + \frac{16\eta_t^2}{n^3}\sum_{i=1}^{n-1} \left \| \sum_{j=i+1}^n \nabla f ( x_* ; \pi^{(t)} ( j ) ) \right \|^2.
\end{align*}

Now taking expectation conditioned on $\mathcal{F}_t$, we get 
\begin{align*}
    \E_t[K_t] 
    &\leq 24L\eta_t^2\E_t \left[F( \tilde{x}_t  ) - F ( x_*)\right]
    + \frac{16\eta_t^2}{n^3}\sum_{i=1}^{n-1} \E_t \left[\left \| \sum_{j=i+1}^n \nabla f ( x_* ; \pi^{(t)} ( j ) ) \right \|^2\right].
\end{align*}
Applying the sample variance Lemma from \citep{mishchenko2020random} we have 
\begin{align*}
    \E_t \left[\left \| \sum_{j=i+1}^n \nabla f ( x_* ; \pi^{(t)} ( j ) ) \right \|^2\right] = \E_t \left[\left \| \sum_{j=i+1}^n \nabla f ( x_* ; \pi^{(t)} ( j ) ) - \nabla F(x_*) \right \|^2\right] \leq \frac{(n-i)i}{n-1} \sigma_{*}^2.
\end{align*}

Substituting this into the previous expression, we get 
\begin{align*}
    \E_t[K_t] 
    &\leq 24L\eta_t^2\E_t \left[F( \tilde{x}_t  ) - F ( x_*)\right]
    + \frac{16\eta_t^2}{n^3}\sum_{i=1}^{n-1} \frac{(n-i)i}{n-1} \sigma_{*}^2\\
    &\leq 24L\eta_t^2\E_t \left[F( \tilde{x}_t  ) - F ( x_*)\right]
    + \frac{16\eta_t^2 \sigma_{*}^2}{3n},
\end{align*}
where we use the facts that $ \sum_{i=1}^{n-1} \frac{i (n-i)}{(n-1)} \leq \frac{n(n+1)}{6} \leq \frac{n^2}{3}$.
Taking total expectation, we have the estimate of Lemma \ref{lem_bound_convex_RR}.
\Eproof

\subsection*{Proof of Theorem \ref{thm_convex_RR}} 

Let us start with inequality \eqref{eq_020_0} from Lemma \ref{lem_bound_common}. Taking total expectation and applying Lemma \ref{lem_bound_convex_RR} we have
\begin{align*}
    T(T+2) \E[F ( \tilde{x}_{T} ) - F ( x_{*} ) ]
    & \leq   \sum_{t=1}^T\frac{L^2\eta_t(t+1)^2}{2 \epsilon_t}  \E[K_t] -\sum_{t=1}^T \E[F ( \tilde{x}_{t} ) - F ( x_{*} ) ]\\
    & \qquad+ \sum_{t=1}^T\frac{2}{\eta_t}  \E \left[\| v^{(t-1)} -  x_{*} \|^2\right] - \sum_{t=1}^T \frac{2}{\eta_t}(1-\epsilon_t) \E \left[\| v^{(t)} -  x_{*} \|^2\right] \\
    & \overset{\eqref{eq_019}}{\leq}    \sum_{t=1}^T\frac{4L^2\eta_t^3(t+1)^2}{\epsilon_t}  \left(3L \E \left[F( \tilde{x}_t  ) - F ( x_*)\right]
    + \frac{2\sigma_{*}^2}{3n} \right) -\sum_{t=1}^T \E[F ( \tilde{x}_{t} ) - F ( x_{*} ) ]\\
    & \qquad+ \sum_{t=1}^T\frac{2}{\eta_t}  \E \left[\| v^{(t-1)} -  x_{*} \|^2\right] - \sum_{t=1}^T \frac{2}{\eta_t}(1-\epsilon_t) \E \left[\| v^{(t)} -  x_{*} \|^2\right] .
\end{align*}
From the choice $\eta_t = \frac{k\alpha^t}{LT}$ we have $\frac{2}{\eta_t} = \frac{2LT}{k\alpha^t}$ and
\begin{align*}
     T(T+2) \E[F ( \tilde{x}_{T} ) - F ( x_{*} ) ]&\leq   \sum_{t=1}^T\frac{k^3\alpha^{3t}}{L^3T^3} \frac{4L^2(t+1)^2}{\epsilon_t}  \left(3L \E \left[F( \tilde{x}_t  ) - F ( x_*)\right]     + \frac{2\sigma_{*}^2}{3n} \right) -\sum_{t=1}^T \E [F ( \tilde{x}_t ) - F ( x_{*} ) ]\\
    & \quad+ \sum_{t=1}^T \frac{2LT}{k\alpha^t} \E \left[ \| v^{(t-1)} -  x_{*} \|^2\right] - \sum_{t=1}^T  \frac{2LT}{k\alpha^t}(1-\epsilon_t) \E \left[\| v^{(t)} -  x_{*} \|^2\right].
\end{align*}
In addition, we choose $\epsilon_t = \frac{\alpha -1}{\alpha}$ and $(1-\epsilon_t) = \frac{1}{\alpha}$. The last two terms cancel out that
\begin{align*}
     T(T+2) \E[F ( \tilde{x}_{T} ) - F ( x_{*} ) ]
    &\leq   \sum_{t=1}^T\frac{k^3\alpha^{3t}}{L^3T^3} \frac{4\alpha L^2(t+1)^2}{\alpha -1}  \left(3L \E \left[F( \tilde{x}_t  ) - F ( x_*)\right]     + \frac{2\sigma_{*}^2}{3n} \right) -\sum_{t=1}^T \E[F ( \tilde{x}_t ) - F ( x_{*} ) ]\\
    & \quad+ \sum_{t=1}^T \frac{2LT}{k\alpha^t}   \E \left[\| v^{(t-1)} -  x_{*} \|^2 \right]- \sum_{t=1}^T  \frac{2LT}{k\alpha^{t+1}}  \E \left[\| v^{(t)} -  x_{*} \|^2\right].\\
    &\leq   \sum_{t=1}^T\frac{k^3\alpha^{3t+1}}{LT^3} \frac{4(t+1)^2}{\alpha-1}  \left(3L \E \left[F( \tilde{x}_t  ) - F ( x_*)\right]     + \frac{2\sigma_{*}^2}{3n}\right) -\sum_{t=1}^T \E[F ( \tilde{x}_t ) - F ( x_{*} ) ]\\
    & \qquad+ \frac{2LT}{k\alpha}  \| v^{0} -  x_{*} \|^2.
\end{align*}
Note that $\alpha = 1+ \frac{1}{T}$ ($1 \leq \alpha \leq \frac{3}{2}$ for $T \geq 2$). Hence $\alpha -1= \frac{1}{T}$, $\alpha^t \leq \alpha^T = \left( 1+ \frac{1}{T}\right)^T \leq e$ and
\begin{align*}
     T(T+2) \E[F ( \tilde{x}_{T} ) - F ( x_{*} ) ]
    &\leq   \sum_{t=1}^T\frac{k^3e^3\alpha}{LT^2}  4 (t+1)^2 \left(3L \E \left[F( \tilde{x}_t  ) - F ( x_*)\right]     + \frac{2\sigma_{*}^2}{3n} \right) -\sum_{t=1}^T \E[F ( \tilde{x}_t ) - F ( x_{*} ) ]\\
    & \qquad
    + \frac{2LT}{k\alpha}  \| v^{0} -  x_{*} \|^2 \\
    &\leq   \sum_{t=1}^T \left[\frac{12k^3e^3\alpha(t+1)^2}{T^2}    -1\right] \E[F ( \tilde{x}_t ) - F ( x_{*} ) ] + \sum_{t=1}^T\frac{8k^3e^3\alpha (t+1)^2  \sigma_*^2 }{3nLT^2}   
    \\ & \qquad + \frac{2LT}{k\alpha}  \| v^{0} -  x_{*} \|^2. 
\end{align*}

From the choice $k = \frac{1}{e \alpha \sqrt[3]{12}}$, we have $12k^3 e^3 \alpha^3 = 1$.
Hence for every $t \geq 1$ we have 
\begin{align*}
    \frac{12k^3e^3\alpha(t+1)^2}{T^2} -1 \leq \frac{12k^3e^3\alpha(T+1)^2}{T^2} -1  \leq 12k^3e^3 \alpha^3 -1 = 0,
\end{align*}
where we use the fact that $\alpha = 1+ \frac{1}{T} =  \frac{T+1}{T}$.

We further have 
\begin{align*}
     T(T+2) \E[F ( \tilde{x}_{T} ) - F ( x_{*} ) ]
    &\leq  \sum_{t=1}^T\frac{8k^3e^3 \alpha  (t+1)^2  \sigma_*^2 }{3nLT^2}   
    + \frac{2LT}{k\alpha}  \| v^{0} -  x_{*} \|^2 \\
    &\leq \frac{8k^3e^3 \alpha  (T+2)^3  \sigma_*^2 }{9nLT^2}   
    + \frac{2LT}{k\alpha}  \| v^{0} -  x_{*} \|^2,
\end{align*}
where we use the fact that $\sum_{t=1}^{T} (t+1)^2 \leq \frac{(T+2)^3}{3}$. Dividing both sides by $T(T+2)$ and substituting $k = \frac{1}{e \alpha \sqrt[3]{12}}$ and $12k^3 e^3 \alpha^3 = 1$ we have
\begin{align*}
    \E[F ( \tilde{x}_{T} ) - F ( x_{*} ) ]
    &\leq \frac{8k^3e^3 \alpha  (T+2)^2  \sigma_*^2 }{9nLT^3}   
    + \frac{2L}{k\alpha(T+2)}  \| v^{0} -  x_{*} \|^2\\
    &\leq \frac{ 2(T+2)^2  \sigma_*^2 }{27n \alpha^2 LT^3}   
    + \frac{2Le \sqrt[3]{12}}{T+2}  \| v^{0} -  x_{*} \|^2\\
    &\leq \frac{8 \sigma_*^2 }{27 n LT}   
    + \frac{2Le \sqrt[3]{12}}{T}  \| v^{0} -  x_{*} \|^2, 
\end{align*}
where $(T+2)^2 \leq 4T^2$ for $T \geq 2$. Note that $v^{0} = \tilde{x}_0$, we get the desired results.
\Eproof

\subsection*{Proof of Corollary \ref{cor_comp_RR}: Computational complexity of Theorem \ref{thm_convex_RR}} 
\begin{corollary}\label{cor_comp_RR}
Assume the same conditions as in Theorem \ref{thm_convex_RR}, i.e. Assumption \ref{ass_smooth} and \ref{ass_convex} holds for \eqref{ERM_problem_01} and a randomized schemes is applied. The computational complexity needed by Algorithm~\ref{shuffling_nesterov_02} to reach an $\epsilon$-accurate solution $x$ that satisfies $\E[F(x) - F(x_*)] \leq \epsilon$ is
\begin{align}
    nT
    = \Ocal \left( \frac{\sigma_*^2}{L\epsilon}+ \frac{nL\| \tilde{x}_0 -  x_{*} \|^2}{\epsilon}\right).
    \label{eq_comp_RR}
\end{align} 
\end{corollary}
By Theorem \ref{thm_convex_RR} we have 
\begin{align*}
    \E[F ( \tilde{x}_{T} ) - F ( x_{*} ) ]
    &\leq \frac{8 \sigma_*^2 }{27n LT}   
    + \frac{2Le \sqrt[3]{12}}{T}  \| \tilde{x}_0 -  x_{*} \|^2.
\end{align*} 
In order to reach an $\epsilon$-accurate solution $x = \tilde{x}_{T}$ that satisfies $\E[F(x) - F(x_*)] \leq \epsilon$, we need
\begin{align*}
    \frac{8 \sigma_*^2 }{27n LT} 
     &\leq \frac{\epsilon}{2} \text{ and  }
    \frac{2Le \sqrt[3]{12}}{T}  \| \tilde{x}_0 -  x_{*} \|^2\leq \frac{\epsilon}{2},
\end{align*} 
which is equivalent to 
\begin{align*}
    T
     &\geq \frac{16\sigma_*^2}{27nL\epsilon} \text{ and  }
    T  \geq \frac{4Le \sqrt[3]{12}\| \tilde{x}_0 -  x_{*} \|^2}{\epsilon}.
\end{align*} 
Hence the number of individual gradient evaluations needed is
\begin{align*}
    nT
     = \max \left( \frac{16\sigma_*^2}{27L\epsilon} , \frac{4nLe \sqrt[3]{12}\| \tilde{x}_0 -  x_{*} \|^2}{\epsilon}\right)\leq  \frac{16\sigma_*^2}{27L\epsilon}+ \frac{4nLe \sqrt[3]{12}\| \tilde{x}_0 -  x_{*} \|^2}{\epsilon}= \Ocal \left( \frac{\sigma_*^2}{L\epsilon}+ \frac{nL\| \tilde{x}_0 -  x_{*} \|^2}{\epsilon}\right).
\end{align*} 
\Eproof

\section{Improved Convergence Rate with Initial Condition}
In this section, we propose an initial condition which requires the iterate of our algorithm to be in a small neighborhood of the optimal point. Let us note that the minimizer of $F$ may not be unique, hence the condition holds for some minimizer $x_*$. 
\begin{assumption}\label{ass_bounded_init_weight}
Let $\tilde{x}_0$ the initial point and $E>0$ be a constant. There exists a minimizer $x_*$ of $F$ which satisfies
\begin{align*}
    \| \tilde{x}_0 -  x_{*} \| \leq \frac{E}{\sqrt{n}}.
\end{align*}
\end{assumption}
Although in practice this assumption can be strong, we believe it provides some theoretical insights to investigate the behaviour of SGD Shuffling-type algorithms when they reach a small neighborhood of the minimizer. Our next two Corollaries demonstrates this fact for unified shuffling and randomized schemes respectively. 

\begin{corollary}\label{cor_unified}
Assume the same conditions as in Theorem~\ref{thm_convex_1}, i.e. Assumption \ref{ass_smooth} and \ref{ass_convex} hold for \eqref{ERM_problem_01}. In addition, we assume Assumption \ref{ass_bounded_init_weight} holds for the initial point $\tilde{x}_0$ of in Algorithm~\ref{shuffling_nesterov_02}. Let $\sets{x_i^{(t)}}$ be generated by  Algorithm~\ref{shuffling_nesterov_02} with parameter $\gamma_t = \frac{t-1}{t+2}$, the learning rate $\eta_i^{(t)} := \frac{\eta_t}{n} > 0$ for $\eta_t = \frac{k\alpha^t}{LT} \leq \frac{1}{L}$ where $k = \frac{1}{e \alpha  n^{1/4} \sqrt[3]{12}} > 0$ and $\alpha = 1+ \frac{1}{T}>0$. Then for $T \geq 2$ we have
\begin{align}
     F ( \tilde{x}_{T} ) - F ( x_{*} ) 
    \leq \frac{4 \sigma_*^2 }{9  Ln^{3/4}T}   
    + \frac{2LE^2e \sqrt[3]{12}}{n^{3/4}T}.
    \label{eq_cor_unified}
\end{align} 
\end{corollary}
The convergence rate of Corollary \ref{cor_unified} is expressed as
\begin{equation*}
\Ocal\left(\frac{ \sigma_*^2/L + LE^2  }{n^{3/4}T} \right),
\end{equation*}
which has an improvement of $n^{3/4}$ over the plain setting of Theorem \ref{thm_convex_1}.

\subsection*{Proof of Corollary \ref{cor_unified}} 

We start with the derivation from Theorem \ref{thm_convex_1}
\begin{align*}
     T(T+2) [F ( \tilde{x}_{T} ) - F ( x_{*} ) ]
    &\leq   \sum_{t=1}^T \left[\frac{12k^3e^3\alpha(t+1)^2}{T^2}    -1\right] [F ( \tilde{x}_{t} ) - F ( x_{*} ) ] + \sum_{t=1}^T\frac{4k^3e^3\alpha (t+1)^2  \sigma_*^2 }{LT^2}   
    + \frac{2LT}{k\alpha}  \| v^{0} -  x_{*} \|^2 .
\end{align*}
Note that $v^{0} = \tilde{x}_0$, by Assumption \ref{ass_bounded_init_weight} we have
\begin{align*}
     T(T+2) [F ( \tilde{x}_{T} ) - F ( x_{*} ) ]
    &\leq   \sum_{t=1}^T \left[\frac{12k^3e^3\alpha(t+1)^2}{T^2}    -1\right] [F ( \tilde{x}_{t} ) - F ( x_{*} ) ] + \sum_{t=1}^T\frac{4k^3e^3\alpha (t+1)^2  \sigma_*^2 }{LT^2}   
    + \frac{2LTE^2}{kn\alpha}  .
\end{align*}
From the choice $k = \frac{1}{e \alpha n^{1/4} \sqrt[3]{12}}$, we have $12k^3 e^3 \alpha^3 = \frac{1}{n^{3/4}} \leq 1$.
Hence for every $t \geq 1$ we have 
\begin{align*}
    \frac{12k^3e^3\alpha(t+1)^2}{T^2} -1 \leq \frac{12k^3e^3\alpha(T+1)^2}{T^2} -1  \leq 12k^3e^3 \alpha^3 -1 = 0,
\end{align*}
where we use the fact that $\alpha = 1+ \frac{1}{T} =  \frac{T+1}{T}$.

We further have 
\begin{align*}
     T(T+2) [F ( \tilde{x}_{T} ) - F ( x_{*} ) ]
    &\leq  \sum_{t=1}^T\frac{4k^3e^3 \alpha  (t+1)^2  \sigma_*^2 }{LT^2}   
    + \frac{2LTE^2}{kn\alpha} \\
    &\leq \frac{4k^3e^3 \alpha  (T+2)^3  \sigma_*^2 }{3LT^2}   
    + \frac{2LTE^2}{kn\alpha},
\end{align*}
where we use the fact that $\sum_{t=1}^{T} (t+1)^2 \leq \frac{(T+2)^3}{3}$. Dividing both sides by $T(T+2)$ and substituting $k = \frac{1}{e \alpha n^{1/4} \sqrt[3]{12}}$ and $12k^3 e^3 \alpha^3 = \frac{1}{n^{3/4}}$ we have
\begin{align*}
    F ( \tilde{x}_{T} ) - F ( x_{*} ) 
    &\leq \frac{4k^3e^3 \alpha  (T+2)^2  \sigma_*^2 }{3LT^3}   
    + \frac{2LTE^2}{kn\alpha}\\
    &\leq \frac{ (T+2)^2  \sigma_*^2 }{9 \alpha^2 Ln^{3/4}T^3}   
    + \frac{2LE^2e \sqrt[3]{12} }{n^{3/4}(T+2)} \\
    &\leq \frac{4 \sigma_*^2 }{9  Ln^{3/4}T}   
    + \frac{2LE^2e \sqrt[3]{12}}{n^{3/4}T} , 
\end{align*}
where $(T+2)^2 \leq 4T^2$ for $T \geq 2$. Thus we get the desired results.
\Eproof
\begin{corollary}\label{cor_RR}
Assume the same conditions and parameter setting as in Theorem \ref{thm_convex_RR}, i.e. Assumption \ref{ass_smooth} and \ref{ass_convex}  hold for \eqref{ERM_problem_01}. In addition, we assume Assumption \ref{ass_bounded_init_weight} holds for the initial point $\tilde{x}_0$ of Algorithm~\ref{shuffling_nesterov_02}. Let $\sets{x_i^{(t)}}$ be generated by  Algorithm~\ref{shuffling_nesterov_02} under a \textbf{randomized scheme} with parameter $\gamma_t = \frac{t-1}{t+2}$, the learning rate $\eta_i^{(t)} := \frac{\eta_t}{n} > 0$ for $\eta_t = \frac{k\alpha^t}{LT} \leq \frac{1}{L}$ where $k = \frac{1}{e \alpha \sqrt[3]{12}} > 0$ and $\alpha = 1+ \frac{1}{T}>0$. Then for $T \geq 2$ we have
\begin{align}
    \E[F ( \tilde{x}_{T} ) - F ( x_{*} ) ]
    &\leq \frac{8 \sigma_*^2 }{27n LT}   
    + \frac{2LE^2e \sqrt[3]{12}}{nT}.  \label{eq_cor_convex_RR}
\end{align}
\end{corollary}
The convergence rate of Corollary \ref{cor_RR} is expressed as
\begin{equation*}
\Ocal\left(\frac{ \sigma_*^2/L + LE^2  }{nT} \right),
\end{equation*}
which shows an improvement of $n$ over the standard setting thanks to the application of Randomized schemes and Assumption \ref{ass_bounded_init_weight}. The proof of Corollary \ref{cor_RR} follows straightforwardly from Theorem \ref{thm_convex_2} and as a result, this bound is in expectation form, which is weaker than the deterministic criteria in Corollary \ref{cor_unified}.

\section{Detailed Implementation and Additional Experiments}\label{sec:experiments_details}
In this section, we explain the detailed hyper-parameter tuning in Section~\ref{sec:experiments}.

\subsection{Experiment Settings}
For the binary classification experiment, we consider the following convex problem:
\begin{align*}
   \min_{w \in \mathbb{R}^d} \Big \{ F(w) \! := \! \frac{1}{n} \sum_{i=1}^n \log(1 \! + \! \exp(- y_i x_i^\top w )) \!  \Big \}, 
\end{align*}
where $\sets{(x_i, y_i)}_{i=1}^n$ is a set of training samples with $x_i \in \mathbb{R}^d$ and $y_i \in \{-1,1\}$. 
For the image classification, we experiment with the following minimization problem:
\begin{align*}
    \min_{w \in \mathbb{R}^d} \Big \{F(w):= - \frac{1}{n}\sum_{i=1}^n y_i^{\top} \log( \text{softmax} ( h(w ; i) ) ) \Big \},
\end{align*}
where $h(\cdot;i)$ can be convex or non-convex. The input data $\sets{x_i}_{i=1}^n$ are in $\R^d$ and the output labels $\sets{y_i}_{i=1}^n$ are one-hot vectors in $\R^c$, where $c$ is the number of classes. 
Note that this problem can be written as  $f(w ; i) = \phi_i ( h ( w  ; i) ) $ where $\phi_i$ is the convex softmax function \citep{Nguyen2022_OptDL}.

\subsection{Comparing NASG with Other Methods}\label{sec:compare1}

For the motivational experiment in Section \ref{sec:shuffling_nesterov}, we use the same setting as the binary classification in Section \ref{subsec:exp_binary}. 

At the tuning stage, we test each method for 20 epochs. 
We run every algorithm with a constant learning rate where the learning rates follows a grid search and select the ones that perform best according to their results. These hyperparameters are choosen for the main training stage that lasts 100 and 200 epochs (for binary experiment and image classification, respectively). 
The hyper-parameters tuning strategy for our main experiments is given below:
\begin{compactitem}
    \item For NASG the searching grid is $\{1, 0.5, 0.1, 0.05, 0.01, 0.005, 0.001\}$.
    \item For  deterministic NAG, the searching grid is $\{50, 10, 5, 1, 0.5, 0.1, 0.05, 0.01, 0.005, 0.001\}$.
    \item For  NASG-PI (applying the Nesterov momentum term for each iteration), the searching grid is $\{10, 5, 1, 0.5, 0.1, 0.05, 0.01, 0.005, 0.001\}$. We describe this method in Algorithm \ref{shuffling_nesterov}. 
\end{compactitem}
\begin{algorithm}[hpt!]
   \caption{Nesterov Accelerated Shuffling Gradient - Per Iteration
(NASG - PI) }\label{shuffling_nesterov}
\begin{algorithmic}[1]
   \STATE {\bfseries Initialization:} Choose an initial point $\tilde{x}_0, \tilde{y}_0 \in \mathbb{R}^d$.
   \FOR{$t=1,2,\cdots,T $}
   \STATE Set $x_0^{(t)} := \tilde{x}_{t-1}$ and $y_0^{(t)} := \tilde{y}_{t-1}$;
   \STATE Generate any permutation  $\pi^{(t)}$ of $[n]$ (either deterministic or random);
   \FOR{$i = 1,\cdots, n$}
    \STATE Update $x_{i}^{(t)} := y_{i-1}^{(t)} - \eta_i^{(t)} \nabla f ( y_{i-1}^{(t)} ; \pi^{(t)} ( i ) )$; 
    \STATE Update $y_{i}^{(t)} := x_{i}^{(t)} + \frac{t-1}{t+2} ( x_{i}^{(t)} - x_{i-1}^{(t)} )$; 
   \ENDFOR
   \STATE Set $\tilde{x}_t := x_{n}^{(t)}$ and $\tilde{y}_t := y_{n}^{(t)}$;
   \ENDFOR
\end{algorithmic}
\end{algorithm}

In the two main sets of algorithm, we compare our algorithm with Stochastic Gradient Descent (SGD) and two other methods: SGD with Momentum (SGD-M) \citep{Polyak1964} and Adam \citep{Kingma2014}. 
To have a fair comparison, a random reshuffling strategy is applied to all methods.

At the tuning stage, we test each method for 20 epochs. 
We run every algorithm with a constant learning rate where the learning rates follows a grid search and select the ones that perform best according to their results. These hyperparameters are choosen for the main training stage that lasts 100 and 200 epochs (for binary experiment and image classification, respectively). 
The hyper-parameters tuning strategy for our main experiments is given below:
\begin{compactitem}
    \item For SGD and NASG the searching grid is $\{1, 0.5, 0.1, 0.05, 0.01, 0.005, 0.001\}$.
    \item For SGD-M, we update the weights using the following rule: 
    \begin{align*}
        &m_{i+1}^{(t)} := \beta m_{i}^{(t)} + g_i^{(t)}\\
        &w_{i+1}^{(t)} := w_{i}^{(t)} - \eta_{i}^{(t)} m_{i+1}^{(t)},
    \end{align*}
    where $g_i^{(t)}$ is the $(i+1)$-th gradient at epoch $t$. 
    Note that this momentum update is implemented in PyTorch with the default value $\beta = 0.9$. 
    Hence, we choose this setting for SGD-M, and we tune the learning rate using the grid search as in the SGD algorithm. 
    \item For Adam, we fixed two hyper-parameters $\beta_1 := 0.9$, $\beta_2 := 0.999$ as in the original paper. 
    Since the default learning rate for Adam is $0.001$, we let our searching grid be $\{0.005, 0.001, 0.0005\}$. 
    We note that since the best learning rate for Adam is usually $0.001$, its hyper-parameter tuning process requires little effort than other algorithms in our experiments. 
\end{compactitem}

\end{document}